\documentclass[a4paper,10pt,leqno]{article}
\usepackage{amsmath, amsfonts, amssymb, amsthm}
\numberwithin{equation}{section}

\begin{document}

\title{ On the Cauchy problem for a two-component
Degasperis-Procesi system}

\author{Kai Yan\footnote{email:
yankai419@163.com }\quad and\quad Zhaoyang Yin\footnote{email: mcsyzy@mail.sysu.edu.cn}\\
 Department of Mathematics,
Sun Yat-sen University,\\510275 Guangzhou, China }
\date{}
\maketitle

\begin{abstract}
This paper is concerned with the Cauchy problem for a two-component
Degasperis-Procesi system. Firstly, the local well-posedness for this system in the nonhomogeneous Besov spaces is established. Then the precise blow-up scenario for strong solutions to the system is derived. Finally, two new blow-up criterions and the exact blow-up rate of strong solutions to the system are presented.\\


\end{abstract}

\section{Introduction}
\par
In this paper we consider the Cauchy problem of the following
two-component  Degasperis-Procesi system:
\begin{equation}
\left\{\begin{array}{ll}m_{t}+3m{u_{x}}+m_{x}u+k_{3}\rho\rho_{x}=0,&t > 0,\,x\in \mathbb{R},\\
 \rho_{t}+k_{2} u\rho_{x}+(k_{1}+k_{2})u_{x}\rho=0, &t > 0,\,x\in \mathbb{R},\\
u(0,x) = u_{0}(x),&x\in \mathbb{R}, \\
\rho(0,x) = \rho_{0}(x),&x\in \mathbb{R},\end{array}\right.
\end{equation}
where $m=u-u_{xx}$, while $(k_{1},k_{2},k_{3})=(1,1,c)$ or $(c,1,0)$ and $c$ takes an arbitrary value.

Using the Green's function $p(x)\triangleq \frac{1}{2}e^{-|x|}$, $x\in \mathbb{R}$ and the identity
$(1-\partial^{2}_{x})^{-1}f = p\ast f $ for all $f \in L^{2}(\mathbb{R})$, we can rewrite System (1.1) as follows:
\begin{equation}
\left\{\begin{array}{ll}u_t+u u_x= P(D)(\frac 3 2 {u^2}+\frac {k_3}{2} {\rho^2}),&t > 0,\,x\in \mathbb{R},\\
 \rho_{t}+k_{2} u\rho_{x}=-(k_{1}+k_{2})u_{x}\rho, &t > 0,\,x\in \mathbb{R},\\
u(0,x) = u_{0}(x),&x\in \mathbb{R}, \\
\rho(0,x) = \rho_{0}(x),&x\in \mathbb{R},\end{array}\right.
\end{equation}
where the operator $P(D)\triangleq -\partial_x (1-\partial^2_x)^{-1}$.

\par
 System (1.1) as the Hamiltonian extension of the Degasperis-Procesi equation
was firstly proposed in \cite{Po1} and it is generated by the Hamiltonian operator
which is a Dirac-reduced operator of the generalized but degenerated second
Hamiltonian operator of the Boussinesq equation. As the author mentioned, he can not verify the
integrability of System (1.1) in his way \cite{Po1}. However, it does not mean that this system is
not integrable. In particular, for $(k_{1},k_{2},k_{3})=(c,1,0)$, System (1.1) is no more coupled and the equation
on $\rho$ becomes linear. Therefore, we only consider the case $(k_{1},k_{2},k_{3})=(1,1,c)$ in the present paper. That is
\begin{equation}
\left\{\begin{array}{ll}u_t+u u_x=P(D)(\frac 3 2 {u^2}+\frac {c}{2} {\rho^2}),&t > 0,\,x\in \mathbb{R},\\
 \rho_{t}+ u\rho_{x}=-2u_{x}\rho, &t > 0,\,x\in \mathbb{R},\\
u(0,x) = u_{0}(x),&x\in \mathbb{R}, \\
\rho(0,x) = \rho_{0}(x),&x\in \mathbb{R}.\end{array}\right.
\end{equation}

For $\rho\equiv 0$, System (1.1) becomes the Degasperis-Procesi
equation \cite{D-P}. It was proved formally integrable by constructing a Lax pair \cite{D-H-H}. Moreover, they also presented \cite{D-H-H} that the DP equation has a bi-Hamiltonian structure and an infinite number of conservation laws,
and admits exact peakon solutions which are analogous to the Camassa-Holm peakons \cite{C-H,C-S,C-S1}. The DP equation can be
viewed as a model for nonlinear shallow water dynamics and its
asymptotic accuracy is the same as for the Camassa-Holm shallow
water equation \cite{C-H,C-L,D-G-H1,D-G-H2,Ip}. It can be obtained from the shallow water elevation equation by an appropriate Kodama transformation \cite{D-G-H2}. The numerical  stability of solitons and peakons, the multi-soliton solutions and their peakon limits, together with an inverse scattering method to compute $n$-peakon solutions to DP equation have been investigated respectively in \cite{H-S,L-S,MY}. Furthermore, the traveling wave solutions and the classification of all weak traveling wave solutions to DP equation were presented in \cite{Le,V-P}.  After the DP equation appeared, it has been studied in many works \cite{C-K,E-Y,H-D,Le,Lu,MY,Y0,Y3}. For example, the author established the local well-posedness to DP equation with initial data $u_0 \in H^s(\mathbb{R}), s>{\frac 3 2}$ on the line \cite{Y3} and on the circle \cite{Y0}, and derived the precise blow-up scenario and a blow-up result. The global existence of strong solutions and global weak solutions to DP equation were shown in \cite{Y4,Y5}.
Similar to the Camassa-Holm equation \cite{Cf,C-Ep,C-E,C-M,Y1,Y2}, the DP
equation has not only global strong solutions \cite{L-Y,Y4} but also
blow-up solutions \cite{E-L-Y,E-L-Y1,L-Y,Y4}. Apart from these, it  has global entropy weak solutions in $L^{1}(\mathbb{R})\cap BV(\mathbb{R})$ and
$L^{2}(\mathbb{R})\cap L^{4}(\mathbb{R})$, cf. \cite{C-K0}.

Although the DP equation is very similar to the
Camass-Holm equation in many aspects, especially in the structure of equation, there are some essential differences between the two equations.
One of the famous features of DP equation is that it has not
only peakon solutions $u_c (t,x)=ce^{-|x-ct|}$ with $c>0$ \cite{D-H-H} and periodic peakon solutions
\cite{Y5}, but also shock peakons \cite{Lu} and the periodic shock waves \cite{E-L-Y1}. Besides, the
Camass-Holm equation is a re-expression of geodesic flow on the diffeomorphism group \cite{C-K} or on the Bott-Virasoro group \cite{Mi}, while the DP equation can be regarded as a non-metric Euler equation \cite{E-B}.

Recently, a large amount of literature was devoted to the two-component Camass-Holm system \cite{C-L-Z,C-I,E-Le-Y,Fa,G-Y,G-Y1,G-L1,G-L2,Iz,Y-Y}. It is noted that the authors in \cite{Y-Y} studied the analytic solutions of the Cauchy problem for two-component Camass-Holm shallow water systems, which were proved in both variables, globally in space and locally in time. The used approach in \cite{Y-Y} depends strongly on the structure of the given system, so we can obtain the same analyticity results of System (1.3) as those in two-component Camass-Holm system. However, the Cauchy problem of System (1.3) in Besov spaces has not been discussed yet. The goal of this paper is to establish the local well-posedness of System (1.3) in the nonhomogeneous Besov spaces, derive the precise blow-up scenario of strong solutions to the system, and give the new blow-up criterions with respect to the initial data and the exact blow-up rate of strong solutions to the system. Most of our results can be carried out to the periodic case and to homogeneous Besov spaces.

To solve the problem, we mainly use the ideas of \cite{D1,G-L1,G-L2,L-Y}. One of the difficulties is the treatment of critical index in proving local well-posedness of System (1.3), which has been overcome by the interpolation method in some sense. On the other hand, the ${H^1}\times {L^2}$-norm conserved quantity plays a key role in studying the blow-up phenomenon of the two-component Camass-Holm system \cite{G-L1,G-L2}. Unfortunately, one can not find this similar conservation law of System (1.3). This difficulty has been dealt with in some sense by obtaining a priori estimate $L^\infty$-norm of the first component of the solutions to System (1.3) and making good use of the structure of the system itself.

We now conclude this introduction by outlining the rest of the paper. In Section 2, we will recall  some facts on the Littlewood-Paley decomposition, the nonhomogeneous Besov spaces and their some useful properties, and the transport equation theory. In Section 3, we establish the local well-posedness of the system. In Section 4, we derive the precise blow-up scenario for strong solutions to the system. Section 5 is devoted to some new blow-up results and the exact blow-up rate of strong solutions to the system.

\section{Preliminaries}
\newtheorem {remark2}{Remark}[section]
\newtheorem{theorem2}{Theorem}[section]
\newtheorem{lemma2}{Lemma}[section]
\newtheorem{definition2}{Definition}[section]
\newtheorem{proposition2}{Proposition}[section]

\par
In this section, we will recall some facts on the Littlewood-Paley decomposition, the nonhomogeneous Besov spaces and their some useful properties, and the transport equation theory, which will be used in the sequel.

\begin{proposition2}
\cite{D3} (Littlewood-Paley decomposition) There exists a couple of smooth functions $(\chi,\varphi)$ valued in $[0,1]$, such that $\chi$ is supported in the ball $B\triangleq \{\xi\in\mathbb{R}^n:|\xi|\leq \frac 4 3\}$, and $\varphi$ is supported in the ring $C\triangleq \{\xi\in\mathbb{R}^n:\frac 3 4\leq|\xi|\leq \frac 8 3\}$. Moreover,
$$\forall\,\ \xi\in\mathbb{R}^n,\,\ \chi(\xi)+{\sum\limits_{q\in\mathbb{N}}\varphi(2^{-q}\xi)}=1,$$
and
$$\textrm{supp}\,\ \varphi(2^{-q}\cdot)\cap \textrm{supp}\,\ \varphi(2^{-q^{'}}\cdot)=\emptyset,\,\ if\,\ |q-q^{'}|\geq 2,$$
$$\textrm{supp}\,\ \chi(\cdot)\cap \textrm{supp}\,\ \varphi(2^{-q}\cdot)=\emptyset,\,\ if\,\ q\geq 1.$$
\end{proposition2}

Then for all $u \in \mathcal{S}^{'}$, we can define the nonhomogeneous dyadic blocks as follows. Let
$$\Delta_q{u}\triangleq 0,\,\ if\,\ q\leq -2,$$
$$\Delta_{-1}{u}\triangleq \chi(D)u=\mathcal{F}^{-1}\chi \mathcal{F}u,$$
$$\Delta_q{u}\triangleq \varphi(2^{-q}D)u=\mathcal{F}^{-1}\varphi(2^{-q}\xi)\mathcal{F}u,\,\ if \,\ q\geq 0.$$
Hence,
$$ u={\sum\limits_{q\in\mathbb{Z}}}\Delta_q{u} \quad  in \,\  \mathcal{S}^{'}(\mathbb{R}^n),$$
where the right-hand side is called the nonhomogeneous Littlewood-Paley decomposition of $u$.

\begin{remark2}
(1) The low frequency cut-off  $S_q$ is defined by
$$S_q{u}\triangleq {\sum\limits_{p=-1}^{q-1}}\Delta_p{u}=\chi(2^{-q}D)u=\mathcal{F}^{-1}\chi(2^{-q}\xi)\mathcal{F}u,
\,\ \forall\,\ q\in\mathbb{N}.$$
(2) The Littlewood-Paley decomposition is quasi-orthogonal in $L^2$ in the following sense:
$$\Delta_p\Delta_q{u}\equiv 0,\quad if \,\ |p-q|\geq2,$$
$$\Delta_q(S_{p-1}u \Delta_p{v})\equiv 0,\quad if \,\ |p-q|\geq5,$$
for all $u,v \in \mathcal{S}^{'}(\mathbb{R}^n)$.\\
(3) Thanks to Young's inequality, we get
$$||\Delta_q{u}||_{L^p},\,\ ||S_q{u}||_{L^p}\leq C||u||_{L^p},\quad \forall\,\ 1\leq p\leq\infty,$$
where $C$ is a positive constant independent of $q$.
\end{remark2}

\begin{definition2}
\cite{D3} (Besov spaces) Let $s\in\mathbb{R}, 1\leq p,r\leq\infty$. The nonhomogeneous Besov space $B^s_{p,r}(\mathbb{R}^n)$ ($B^s_{p,r}$ for short) is defined by
$$B^s_{p,r}(\mathbb{R}^n)\triangleq \{f \in \mathcal{S}^{'}(\mathbb{R}^n):||f||_{B^s_{p,r}}< \infty\},$$
where
$$||f||_{B^s_{p,r}}\triangleq ||2^{qs}\Delta_q{f}||_{l^r(L^p)}=||(2^{qs}||\Delta_q{f}||_{L^p})_{q\geq-1}||_{l^r}.$$
If $s=\infty$, $B^\infty_{p,r}\triangleq \bigcap\limits_{s\in\mathbb{R}}B^s_{p,r}$.
\end{definition2}

\begin{definition2}
Let $T>0$, $s\in\mathbb{R}$ and $1\leq p \leq\infty$. Set
$$E^s_{p,r}(T)\triangleq C([0,T]; B^s_{p,r})\cap C^1{([0,T]; B^{s-1}_{p,r})},\quad \text{if}\,\ r<\infty,$$
$$E^s_{p,\infty}(T)\triangleq L^{\infty}([0,T]; B^s_{p,\infty})\cap Lip\,([0,T]; B^{s-1}_{p,\infty})$$
and
$$E^s_{p,r}\triangleq \bigcap\limits_{T>0}E^s_{p,r}(T).$$
\end{definition2}

\begin{remark2}
By Definition 2.1 and Remark 2.1(3), we can deduce that
$$||\Delta_q{u}||_{B^s_{p,r}},\,\ ||S_q{u}||_{B^s_{p,r}}\leq C||u||_{B^s_{p,r}},$$
where $C$ is a positive constant independent of $q$.
\end{remark2}

In the following proposition, we list some important properties of Besov spaces.
\begin{proposition2}
\cite{Cp,D3,G-L2} Suppose that $s\in\mathbb{R}, 1\leq p,r,p_{i},r_{i}\leq\infty, i=1,2.$  We have\\
(1) Topological properties: $B^s_{p,r}$ is a Banach space which is continuously embedded in $\mathcal{S}^{'}$.\\
(2) Density: $\mathcal{C}^{\infty}_c$ is dense in $B^s_{p,r}$ $\Longleftrightarrow 1\leq p,r < \infty.$\\
(3) Embedding: $B^s_{p_1,r_1}\hookrightarrow B^{s-n({1\over p_1}-{1\over p_2})}_{p_2,r_2}$, \,\  if\quad $p_1\leq p_2$ and  $r_1\leq r_2$,
$$B^{s_2}_{p,r_2}\hookrightarrow B^{s_1}_{p,r_1}\quad  locally \,\ compact,\quad if\,\ \, s_1 < s_2.$$
(4) Algebraic properties: $\forall s>0$, $B^s_{p,r}\bigcap L^\infty$ is an algebra.
Moreover, $B^s_{p,r}$ is an algebra, provided that $s>{n\over p}$ or $s\geq{n\over p}\,\ and\,\ r=1$.\\
(5) 1-D Morse-type estimates:\\
(i) For $s>0$,
\begin{eqnarray}||fg||_{B^s_{p,r}(\mathbb{R})}\leq C(||f||_{B^s_{p,r}(\mathbb{R})}||g||_{L^\infty(\mathbb{R})}+||g||_{B^s_{p,r}(\mathbb{R})}||f||_{L^\infty(\mathbb{R})}).
\end{eqnarray}
(ii) $\forall s_{1}\leq \frac{1}{p} < s_2$ ($s_{2}\geq {1\over p}$ if $r=1$) and $s_1+s_2>0$, we have
\begin{eqnarray}
||fg||_{B^{s_1}_{p,r}(\mathbb{R})}\leq C||f||_{B^{s_1}_{p,r}(\mathbb{R})}||g||_{B^{s_2}_{p,r}(\mathbb{R})}.
\end{eqnarray}
(iii) In Sobolev spaces $H^s(\mathbb{R})= B^s_{2,2}(\mathbb{R})$, we have for $s>0$,
\begin{eqnarray}||f\partial_x g||_{H^s(\mathbb{R})}\leq
C(||f||_{H^{s+1}(\mathbb{R})}||g||_{L^\infty(\mathbb{R})}+||f||_{L^\infty(\mathbb{R})} ||\partial_x g||_{H^s(\mathbb{R})}),
\end{eqnarray}
where $C$ is a positive constant independent of $f$ and $g$.\\
(6) Complex interpolation:
\begin{eqnarray}
\quad \quad \quad ||f||_{B^{\theta {s_1}+(1-\theta){s_2}}_{p,r}}\leq ||f||^{\theta}_{B^{s_1}_{p,r}}||f||^{1-\theta}_{B^{s_2}_{p,r}}, \quad \forall u\in B^{s_1}_{p,r}\cap B^{s_1}_{p,r},\quad \forall \theta \in[0,1].
\end{eqnarray}
(7) Fatou lemma: if $(u_n)_{n\in \mathbb{N}}$ is bounded in $B^s_{p,r}$ and $u_n \to u $ in $\mathcal{S}^{'}$, then $u\in B^s_{p,r}$ and
$$||u||_{B^s_{p,r}}\leq \liminf\limits_{n\to \infty} ||u_n||_{B^s_{p,r}}.$$
(8) Let $m\in\mathbb{R}$ and $f$ be a $S^m$-multiplier (i.e., $f: \mathbb{R}^n \to \mathbb{R}$ is smooth and satisfies that $\forall\,\ \alpha\in \mathbb{N}^n,\,\ \exists\,\ a \,\ constant\,\ C_\alpha,\,\ s.t. \,\ |\partial^\alpha{f(\xi)}|\leq C_\alpha(1+|\xi|)^{m-|\alpha|}\,\ for \,\ all\,\ \xi \in \mathbb{R}^n$). Then the operator $f(D)$ is continuous from $B^s_{p,r}$ to $B^{s-m}_{p,r}.$
\end{proposition2}

Now we state some useful results in the transport equation theory, which are crucial to the proofs of our main theorems later.
\begin{lemma2}
\cite{D1,D3} (A priori estimates in Besov spaces) Let $1\leq p,r\leq \infty$ and $s>-\min ({1\over p}, 1-{1\over p}).$ Assume that $f_0\in B^s_{p,r}$, $F\in L^1(0,T; B^s_{p,r})$, and $\partial_x v$ belongs to $L^1(0,T; B^{s-1}_{p,r})$ if $s> 1+{1\over p}$ or to $L^1(0,T; B^{1\over p}_{p,r}\cap L^\infty)$ otherwise. If $f\in L^\infty(0,T; B^s_{p,r})\bigcap C([0,T]; \mathcal{S}^{'})$ solves the following 1-D linear transport equation:

\[(T)\left\{\begin{array}{l}
\partial_t f+v\,\partial_x f=F,\\
f|_{t=0} =f_0.
\end{array}\right.\]
Then there exists a constant $C$ depending only on $s,p$ and $r$, and such that the following statements hold:\\
(1) If $r=1$ or $s\neq 1+{1\over p}$,
\begin{equation*}
||f(t)||_{B^s_{p,r}}\leq ||f_0||_{B^s_{p,r}}\,+\, \int_0^t ||F(\tau)||_{B^s_{p,r}}d\tau\,+\, C\int_0^t V^{'}(\tau)||f(\tau)||_{B^s_{p,r}} d\tau
\end{equation*}
or hence,
\begin{equation*}
||f(t)||_{B^s_{p,r}}\leq e^{CV(t)} (||f_0||_{B^s_{p,r}}\,+\, \int_0^t e^{-CV(\tau)} ||F(\tau)||_{B^s_{p,r}}d\tau)
\end{equation*}
with $V(t)=\int_0^t ||\partial_x v(\tau)||_{B^{1\over p}_{p,r}\cap L^\infty}d\tau$ if $s<1+{1\over p}$ and $V(t)=\int_0^t ||\partial_x v(\tau)||_{B^{s-1}_{p,r}}d\tau$ else.\\
(2) If $s\leq 1+{1\over p}$, and $f^{'}_0\in L^\infty$, $\partial_x f\in L^\infty((0,T])\times \mathbb{R})$ and $\partial_x F\in L^1(0,T;L^\infty)$, then
\begin{eqnarray*}
&& \ \ \ \
||f(t)||_{B^s_{p,r}}+||\partial_x f(t)||_{L^\infty}\\
\nonumber&\leq& e^{CV(t)} (||f_0||_{B^s_{p,r}}\,+\,||\partial_x f_0||_{L^\infty}\,+\, \int_0^t e^{-CV(\tau)} (||F(\tau)||_{B^s_{p,r}}+||\partial_x F(\tau)||_{L^\infty})d\tau),
\end{eqnarray*}
with $V(t)=\int_0^t ||\partial_x v(\tau)||_{B^{1\over p}_{p,r}\cap L^\infty}d\tau$.\\
(3) If $f=v$, then for all $s>0$, (2.2) holds true with $V(t)=\int_0^t ||\partial_x v(\tau)||_{L^\infty}d\tau $.\\
(4) If $r<\infty$, then $f\in C([0,T]; B^s_{p,r})$. If $r=\infty$, then $f\in C([0,T]; B^{s^{'}}_{p,1})$ for all $s^{'}<s$.
\end{lemma2}

\begin{lemma2}
\cite{D3} (Existence and uniqueness) Let $p,r,s,f_0$ and $F$ be as in the statement of Lemma 2.1. Assume that $v\in L^\rho(0,T; B^{-M}_{\infty,\infty})$ for some $\rho >1$ and $M>0$, and $\partial_x v \in L^1(0,T; {B^{s-1}_{p,r}})$ if $s> 1+{1\over p}$ or $s=1+{1\over p}$ and $r=1$, and  $\partial_x v \in L^1(0,T; B^{1\over p}_{p,\infty}\cap L^\infty)$ if $s<1+{1\over p}$.
Then (T) has a unique solution $f\in L^{\infty}(0,T; B^s_{p,r})\bigcap\,\big(\bigcap\limits_{s^{'}<s} C([0,T]; B^{s^{'}}_{p,1})\big)$ and the inequalities of Lemma 2.1 can hold true. Moreover, if $r<\infty$, then $f\in C([0,T]; B^s_{p,r})$.
\end{lemma2}

\begin{lemma2}
\cite{G-L2} (A priori estimate in Sobolev spaces) Let $0<\sigma<1$. Assume that $f_0\in H^\sigma$, $F\in L^1(0,T; H^\sigma)$, and $v,\partial_x v\in L^1(0,T; L^\infty)$. If $f\in L^\infty(0,T; H^\sigma)\bigcap C([0,T]; \mathcal{S}^{'})$ solves $(T)$,
then $f\in C([0,T]; H^{\sigma})$, and there exists a constant $C$ depending only on $\sigma$ such that the following statement holds:
\begin{equation*}
||f(t)||_{H^{\sigma}}\leq ||f_0||_{H^{\sigma}}\,+\, \ C \int_0^t ||F(\tau)||_{H^{\sigma}}d\tau\,+\, C\int_0^t V^{'}(\tau)||f(\tau)||_{H^{\sigma}} d\tau
\end{equation*}
or hence,
\begin{equation*}
||f(t)||_{H^{\sigma}}\leq e^{CV(t)} (||f_0||_{H^{\sigma}}\,+\, \int_0^t ||F(\tau)||_{H^{\sigma}}d\tau)
\end{equation*}
with $V(t)=\int_0^t (||v(\tau)||_{L^\infty}+||\partial_x v(\tau)||_{L^\infty})d\tau$.
\end{lemma2}

\section{Local well-posedness}
\newtheorem {remark3}{Remark}[section]
\newtheorem{theorem3}{Theorem}[section]
\newtheorem{lemma3}{Lemma}[section]
\newtheorem{corollary3}{Corollary}[section]
\newtheorem{proposition3}{Proposition}[section]
\par
In this section, we will establish the local well-posedness of System (1.3) in the nonhomogeneous Besov spaces.

For completeness, we firstly apply the classical Kato's semigroup theory \cite{K1} to obtain the local well-posedness of System (1.3) in Sobolev spaces.
More precisely, we have
\begin{theorem3}
Suppose that $z_{0}\triangleq \left(\begin{array}{c}
                                u_{0} \\
                                \rho_{0} \\
                              \end{array}
                            \right)\in H^{s}(\mathbb{R})\times
H^{s-1}(\mathbb{R}),\,\ and \,\ s\geq 2$. There exists a maximal existence time $T=T(\parallel
z_{0}\parallel_{H^{s}(\mathbb{R})\times H^{s-1}(\mathbb{R})})>0$, and
a unique solution $z\triangleq \left(\begin{array}{c}
                                u \\
                                \rho \\
                              \end{array}
                            \right)$ to
 System (1.3) such that
$$
z=z(\cdot,z_{0})\in C([0,T); H^{s}(\mathbb{R})\times
H^{s-1}(\mathbb{R}))\cap C^{1}([0,T);H^{s-1}(\mathbb{R})\times
H^{s-2}(\mathbb{R})).
$$
Moreover, the solution depends continuously on the initial data,
that is, the mapping $z_0\mapsto z(\cdot,z_0):$
$$ H^{s}(\mathbb{R})\!\times\!
H^{s-1}(\mathbb{R})\!\rightarrow\! C([0,T);
H^{s}(\mathbb{R})\!\times\! H^{s-1}(\mathbb{R}))\cap
C^{1}([0,T);H^{s-1}(\mathbb{R})\!\times\! H^{s-2}(\mathbb{R}))
$$
is continuous.
\end{theorem3}

\begin{proof}
The proof is very similar to that in \cite{E-Le-Y}, so we omit it here.
\end{proof}

\par
Now we pay attention to the case in the nonhomogeneous Besov spaces. Uniqueness and continuity with respect to the initial data in some sense can be obtained by the following a priori estimates.
\begin{lemma3}
Let $1\leq p,r\leq \infty$ and $s>\max (2-\frac{1}{p},1+\frac{1}{p},\frac 3 2)$. Suppose that we are given
$\left(\begin{array}{c}
                                u^{i} \\
                                \rho^{i} \\
                              \end{array}
                            \right)\in L^{\infty}(0,T; B^s_{p,r})\cap C([0,T];\mathcal{S}^{'})\times L^{\infty}(0,T; B^{s-1}_{p,r})\cap C([0,T];\mathcal{S}^{'})$ ($i=1,2$) two solutions of System (1.3) with the initial data
$\left(\begin{array}{c}
                                u^{i}_0 \\
                                \rho^{i}_0 \\
                              \end{array}
                            \right)\in B^s_{p,r}\times B^{s-1}_{p,r}$ ($i=1,2$)
                            and let $u^{12}\triangleq u^2-u^1$ and $\rho^{12}\triangleq \rho^2-\rho^1$.
Then for all $t\in[0,T]$, we have\\
(1) if $s>\max (2-\frac{1}{p},1+\frac{1}{p},\frac 3 2)$, but $s\neq 2+\frac{1}{p}, 3+\frac{1}{p}$, then
\begin{eqnarray}
&&||u^{12}(t)||_{B^{s-1}_{p,r}}+||\rho^{12}(t)||_{B^{s-2}_{p,r}}\\
\nonumber&\leq&(||u^{12}_0||_{B^{s-1}_{p,r}}+||\rho^{12}_0||_{B^{s-2}_{p,r}})\\
\nonumber&&\times e^{C\int_0^t (||u^1(\tau)||_{B^s_{p,r}}+||u^2(\tau)||_{B^s_{p,r}}+||\rho^1(\tau)||_{B^{s-1}_{p,r}}+||\rho^2(\tau)||_{B^{s-1}_{p,r}})d\tau};
\end{eqnarray}
(2) if $s= 2+\frac{1}{p}$, then
\begin{eqnarray*}
&&||u^{12}(t)||_{B^{s-1}_{p,r}}+||\rho^{12}(t)||_{B^{s-2}_{p,r}}\\
\nonumber&\leq& C (||u^{12}_0||_{B^{s-1}_{p,r}}+||\rho^{12}_0||_{B^{s-2}_{p,r}})^{\theta}\times (||u^{1}(t)||_{B^{s}_{p,r}}+||u^{2}(t)||_{B^{s}_{p,r}})^{1-\theta}\\
\nonumber&&\times e^{\theta C\int_0^t (||u^1(\tau)||_{B^s_{p,r}}+||u^2(\tau)||_{B^s_{p,r}}+||\rho^1(\tau)||_{B^{s-1}_{p,r}}+||\rho^2(\tau)||_{B^{s-1}_{p,r}})d\tau}\\
\nonumber&&+(||u^{12}_0||_{B^{s-1}_{p,r}}+||\rho^{12}_0||_{B^{s-2}_{p,r}})\\
\nonumber&&\times e^{C\int_0^t (||u^1(\tau)||_{B^s_{p,r}}+||u^2(\tau)||_{B^s_{p,r}}+||\rho^1(\tau)||_{B^{s-1}_{p,r}}+||\rho^2(\tau)||_{B^{s-1}_{p,r}})d\tau};
\end{eqnarray*}
(3) if $s= 3+\frac{1}{p}$, then
\begin{eqnarray*}
&&||u^{12}(t)||_{B^{s-1}_{p,r}}+||\rho^{12}(t)||_{B^{s-2}_{p,r}}\\
\nonumber&\leq& C (||u^{12}_0||_{B^{s-1}_{p,r}}+||\rho^{12}_0||_{B^{s-2}_{p,r}})^{\theta}\times (||\rho^{1}(t)||_{B^{s-1}_{p,r}}+||\rho^{2}(t)||_{B^{s-1}_{p,r}})^{1-\theta}\\
\nonumber&&\times e^{\theta C\int_0^t (||u^1(\tau)||_{B^s_{p,r}}+||u^2(\tau)||_{B^s_{p,r}}+||\rho^1(\tau)||_{B^{s-1}_{p,r}}+||\rho^2(\tau)||_{B^{s-1}_{p,r}})d\tau}\\
\nonumber&&+(||u^{12}_0||_{B^{s-1}_{p,r}}+||\rho^{12}_0||_{B^{s-2}_{p,r}})\\
\nonumber&&\times e^{C\int_0^t (||u^1(\tau)||_{B^s_{p,r}}+||u^2(\tau)||_{B^s_{p,r}}+||\rho^1(\tau)||_{B^{s-1}_{p,r}}+||\rho^2(\tau)||_{B^{s-1}_{p,r}})d\tau},
\end{eqnarray*}
where $\theta\in(0,1)$.
\end{lemma3}

\begin{proof}
It is obvious that $\left(\begin{array}{c}
                                u^{12} \\
                                \rho^{12} \\
                              \end{array}
                            \right)\in L^{\infty}(0,T; B^s_{p,r})\cap C([0,T];\mathcal{S}^{'})\times L^{\infty}(0,T; B^{s-1}_{p,r})\cap C([0,T];\mathcal{S}^{'})$
solves the following Cauchy problem of the transport equations:
\begin{equation}
\left\{\begin{array}{ll}
\partial_t u^{12}+u^{1}\partial_x u^{12}=F(t,x),\\
\partial_t \rho^{12}+u^{1}\partial_x \rho^{12}=G(t,x),\\
u^{12}|_{t=0}=u^{12}_0\triangleq u^{2}_0-u^{1}_0,\\
\rho^{12}|_{t=0}=\rho^{12}_0\triangleq \rho^{2}_0-\rho^{1}_0,
\end{array}\right.
\end{equation}
where $F(t,x)\triangleq -u^{12}\partial_x u^{2}+P(D)\big(\frac {3}{2} u^{12}(u^1+u^2)+\frac {c}{2}\rho^{12}(\rho^1+\rho^2)\big)$ and
$G(t,x)\triangleq -u^{12}\partial_x \rho^{2}-2(\rho^{12}\partial_x u^{1}+\rho^{2}\partial_x u^{12})$.

Claim. For all $s>\max (1+\frac{1}{p},\frac 3 2)$ and $t\in [0,T]$, we have
\begin{eqnarray*}
&&||F(t)||_{B^{s-1}_{p,r}},\,\ ||G(t)||_{B^{s-2}_{p,r}}\\
\nonumber&\leq&C(||u^{12}(t)||_{B^{s-1}_{p,r}}+||\rho^{12}(t)||_{^{s-2}_{p,r}})\\
\nonumber&&\times (||u^1(t)||_{B^s_{p,r}}+||u^2(t)||_{B^s_{p,r}}+||\rho^1(t)||_{B^{s-1}_{p,r}}+||\rho^2(t)||_{B^{s-1}_{p,r}}),
\end{eqnarray*}
where $C=C(s,p,r,c)$ is a positive constant.\\
Indeed, for $s>1+\frac{1}{p}$, $B^{s-1}_{p,r}$ is an algebra, by Proposition 2.2 (4), we have
$$||u^{12}\partial_x u^{2}||_{B^{s-1}_{p,r}}\leq C||u^{12}||_{B^{s-1}_{p,r}}||\partial_x u^{2}||_{B^{s-1}_{p,r}}\leq C||u^{12}||_{B^{s-1}_{p,r}}||u^2||_{B^s_{p,r}}.$$
Note that $P(D)\in Op(S^{-1})$. According to Proposition 2.2 (8) and (2.2), we obtain
$$||P(D)(\frac {3}{2} u^{12}(u^1+u^2))||_{B^{s-1}_{p,r}}\leq C ||u^{12}||_{B^{s-1}_{p,r}}\big(||u^1||_{B^{s-1}_{p,r}}+||u^2||_{B^{s-1}_{p,r}}\big)$$
and
$$||P(D)(\frac {c}{2}\rho^{12}(\rho^1+\rho^2))||_{B^{s-1}_{p,r}}\leq C ||\rho^{12}||_{B^{s-2}_{p,r}}\big(||\rho^1||_{B^{s-1}_{p,r}}+||\rho^2||_{B^{s-1}_{p,r}}\big),$$
if $\max (1+\frac{1}{p},\frac 3 2)<s\leq {2+\frac{1}{p}}$.\\
Otherwise, these inequalities can also hold true in view of the fact $B^{s-2}_{p,r}$ is an algebra as $s>2+\frac{1}{p}$. Therefore,
\begin{eqnarray*}
||F(t)||_{B^{s-1}_{p,r}}
&\leq& C(||u^{12}(t)||_{B^{s-1}_{p,r}}+||\rho^{12}(t)||_{^{s-2}_{p,r}})\\
&\times& (||u^1(t)||_{B^s_{p,r}}+||u^2(t)||_{B^s_{p,r}}+||\rho^1(t)||_{B^{s-1}_{p,r}}+||\rho^2(t)||_{B^{s-1}_{p,r}}),
\end{eqnarray*}
On the other hand, thanks to (2.2), we get
$$||u^{12}\partial_x \rho^{2}||_{B^{s-2}_{p,r}}\leq C ||u^{12}||_{B^{s-1}_{p,r}}||\partial_x \rho^{2}||_{B^{s-2}_{p,r}},$$
$$||\rho^{12}\partial_x u^{1}||_{B^{s-2}_{p,r}}\leq C||\rho^{12}||_{B^{s-2}_{p,r}}||\partial_x u^{1}||_{B^{s-1}_{p,r}} $$
and
$$||\rho^{2}\partial_x u^{12}||_{B^{s-2}_{p,r}}\leq C||\rho^{2}||_{B^{s-1}_{p,r}}||\partial_x u^{12}||_{B^{s-2}_{p,r}},$$
if $\max (1+\frac{1}{p},\frac 3 2)<s\leq {2+\frac{1}{p}}$.\\
For $s>2+\frac{1}{p}$, we can handle it in a similar way. Therefore,
\begin{eqnarray*}
||G(t)||_{B^{s-2}_{p,r}}
&\leq& C(||u^{12}(t)||_{B^{s-1}_{p,r}}+||\rho^{12}(t)||_{^{s-2}_{p,r}})\\
&\times& (||u^1(t)||_{B^s_{p,r}}+||u^2(t)||_{B^s_{p,r}}+||\rho^1(t)||_{B^{s-1}_{p,r}}+||\rho^2(t)||_{B^{s-1}_{p,r}}).
\end{eqnarray*}
This proves our Claim.

Applying Lemma 2.1 (1) and the fact that $||\partial_x w(t)||_{B^{1\over p}_{p,r}\cap L^\infty}\leq C ||w(t)||_{B^s_{p,r}}$,
$||\partial_x w(t)||_{B^{s-3}_{p,r}}\leq C ||\partial_x w(t)||_{B^{s-2}_{p,r}}\leq C ||w(t)||_{B^s_{p,r}}$, if $w\in B^s_{p,r}$ with
$s>\max (2-\frac{1}{p},1+\frac{1}{p},\frac 3 2)$, we can obtain, for case (1),
$$||u^{12}(t)||_{B^{s-1}_{p,r}}\leq ||u^{12}_0||_{B^{s-1}_{p,r}}+\int_0^t ||F(\tau)||_{B^{s-1}_{p,r}}d\tau +C \int_0^t ||u^1(\tau)||_{B^s_{p,r}}||u^{12}(\tau)||_{B^{s-1}_{p,r}}d\tau$$
and
$$||\rho^{12}(t)||_{B^{s-2}_{p,r}}\leq ||\rho^{12}_0||_{B^{s-2}_{p,r}}+\int_0^t ||G(\tau)||_{B^{s-2}_{p,r}}d\tau +C \int_0^t ||u^1(\tau)||_{B^s_{p,r}}||\rho^{12}(\tau)||_{B^{s-2}_{p,r}}d\tau,$$
which together with the Claim yield
\begin{eqnarray*}
&&||u^{12}(t)||_{B^{s-1}_{p,r}}+||\rho^{12}(t)||_{B^{s-2}_{p,r}}\\
\nonumber&\leq&||u^{12}_0||_{B^{s-1}_{p,r}}+||\rho^{12}_0||_{B^{s-2}_{p,r}}+C \int_0^t (||u^{12}(\tau)||_{B^{s-1}_{p,r}}+||\rho^{12}(\tau)||_{^{s-2}_{p,r}})\\
\nonumber&&\times (||u^1(\tau)||_{B^s_{p,r}}+||u^2(\tau)||_{B^s_{p,r}}+||\rho^1(\tau)||_{B^{s-1}_{p,r}}+||\rho^2(\tau)||_{B^{s-1}_{p,r}})d\tau.
\end{eqnarray*}
Taking advantage of Gronwall's inequality, we get (3.1).

For the critical case (2) $s= 2+\frac{1}{p}$, we here use the interpolation method to deal with it. Indeed, if we choose
$s_1\in(\max(2-\frac{1}{p},1+\frac{1}{p},\frac 3 2)-1,s-1)$, $s_2\in (s-1,s)$ and $\theta=\frac{s_2-(s-1)}{s_2-s_1} \in(0,1)$,
then $s-1=\theta s_1+(1-\theta s_2)$. According to Proposition 2.2 (6) and the consequence of case (1), we have
\begin{eqnarray*}
&&||u^{12}(t)||_{B^{s-1}_{p,r}}\\
\nonumber&\leq&||u^{12}(t)||^{\theta}_{B^{s_1}_{p,r}}||u^{12}(t)||^{1-\theta}_{B^{s_2}_{p,r}}\\
\nonumber&\leq&(||u^{1}(t)||_{B^{s_2}_{p,r}}+||u^{2}(t)||_{B^{s_2}_{p,r}})^{1-\theta}
(||u^{12}_0||_{B^{s_1}_{p,r}}+||\rho^{12}_0||_{B^{s_1-1}_{p,r}})^\theta\\
\nonumber&&\times e^{\theta C\int_0^t (||u^1(\tau)||_{B^{s_1+1}_{p,r}}+||u^2(\tau)||_{B^{s_1+1}_{p,r}}+||\rho^1(\tau)||_{B^{s_1}_{p,r}}+||\rho^2(\tau)||_{B^{s_1}_{p,r}})d\tau}\\
\nonumber&\leq& C(||u^{12}_0||_{B^{s-1}_{p,r}}+||\rho^{12}_0||_{B^{s-2}_{p,r}})^{\theta}
(||u^{1}(t)||_{B^{s}_{p,r}}+||u^{2}(t)||_{B^{s}_{p,r}})^{1-\theta}\\
\nonumber&&\times e^{\theta C\int_0^t (||u^1(\tau)||_{B^s_{p,r}}+||u^2(\tau)||_{B^s_{p,r}}+||\rho^1(\tau)||_{B^{s-1}_{p,r}}+||\rho^2(\tau)||_{B^{s-1}_{p,r}})d\tau}.
\end{eqnarray*}
On the other hand, thanks to $s-2=\frac{1}{p}<1+\frac{1}{p}$, then the estimate for $||v^{12}(t)||_{B^{s-2}_{p,r}}$ in case (1) can also hold true. Hence, we can get the desired result.

For the critical case (3) $s= 3+\frac{1}{p}$, its proof is very similar to that of case (2). Therefore, we complete our proof of Lemma 3.1.
\end{proof}

We next construct the approximation solutions to System (1.3) as follows.
\begin{lemma3}
 Let $p$ and $r$ be as in the statement of Lemma 3.1. Assume that $s>\max (2-\frac{1}{p},1+\frac{1}{p},\frac 3 2)$ and $s\neq2+\frac{1}{p}$,
 $z_{0}\triangleq
                            \left(\begin{array}{c}
                                u_0 \\
                                \rho_0 \\
                              \end{array}
                            \right)\in B^s_{p,r}\times B^{s-1}_{p,r}$ and
$z^{0}\triangleq
                            \left(\begin{array}{c}
                                u^0 \\
                                \rho^0 \\
                              \end{array}
                            \right)
                            =
\left(\begin{array}{c}
                                0 \\
                                0 \\
                              \end{array}
                            \right)$. Then\\
(1) there exists a sequence of smooth functions $(z^n)_{n\in \mathbb{N}}\triangleq
\left(\begin{array}{c}
                                u^n \\
                                \rho^n \\
                              \end{array}
                            \right)_{n\in \mathbb{N}}$ belonging to $(C(\mathbb{R}^{+}; B^{\infty}_{p,r}))^2$
and solving the following linear transport equations by induction:
\[(T_n)\left\{\begin{array}{l}
\partial_t u^{n+1}+u^n\,\partial_x u^{n+1}=P(D)(\frac 3 2 (u^n)^2+\frac{c}{2}(\rho^n)^2),\\
\partial_t \rho^{n+1}+u^n\,\partial_x \rho^{n+1}=-2\rho^n \partial_x u^n,\\
u^{n+1}|_{t=0}\triangleq u^{n+1}_0(x)=S_{n+1}u_0,\\
\rho^{n+1}|_{t=0}\triangleq \rho^{n+1}_0(x)=S_{n+1}\rho_0.
\end{array}\right.\]
(2) there exists $T>0$ such that the solutions $(z^n)_{n\in \mathbb{N}}$
is uniformly bounded in $E^s_{p,r}(T)\times E^{s-1}_{p,r}(T)$ and a Cauchy sequence in $C([0,T]; B^{s-1}_{p,r})\times C([0,T]; B^{s-2}_{p,r})$, whence it converges to some limit  $z\triangleq
\left(\begin{array}{c}
                                u \\
                                \rho \\
                              \end{array}
                            \right)$ $\in C([0,T]; B^{s-1}_{p,r})\times C([0,T]; B^{s-2}_{p,r})$.
\end{lemma3}

\begin{proof}
Since all the data $S_{n+1}u_0,\,S_{n+1}\rho_0\in B^{\infty}_{p,r}$, it then follows from Lemma 2.2 and by induction with respect to the index $n$ that (1) holds.

To prove (2), applying Remark (2.2) and  simulating  the proof of Lemma 3.1 (1), we obtain that for $s>\max (2-\frac{1}{p},1+\frac{1}{p},\frac 3 2)$ and $s\neq2+\frac{1}{p}$,
\begin{eqnarray}
a_{n+1}(t)\leq C e^{C U^n(t)}\big(A+\int_0^t e^{-C U^n(\tau)}a^2_n(\tau)d\tau\big),
\end{eqnarray}
where $a_n(t)\triangleq ||u^n(t)||_{B^s_{p,r}}+||\rho^n(t)||_{B^{s-1}_{p,r}}$, $A\triangleq ||u_0||_{B^s_{p,r}}+||\rho_0||_{B^{s-1}_{p,r}}$ and $U^n(t)\triangleq \int_0^t ||u^n(\tau)||_{B^s_{p,r}}d\tau$.\\
Choose $0<\, T\, < \frac{1}{2C^2 A}$ and suppose that
\begin{eqnarray}
a_n(t)\leq \frac{CA}{1-2C^2 At},\quad \forall t\in[0,T].
\end{eqnarray}
Noting that $e^{C (U^n(t)-U^n(\tau))}\leq \sqrt{\frac{1-2C^2 A\tau}{1-2C^2 At}}$ and substituting (3.4) into (3.3) yields
\begin{eqnarray*}
&&a_{n+1}(t)\\
\nonumber&\leq&\frac{CA}{\sqrt{1-2C^2 At}}+\frac{C}{\sqrt{1-2C^2 At}} \int_0^t \frac{C^2 A^2}{(1-2C^2 A\tau)^\frac{3}{2}} d\tau\\
\nonumber&=&\frac{CA}{\sqrt{1-2C^2 At}}+\frac{C}{\sqrt{1-2C^2 At}} (\frac{A}{\sqrt{1-2C^2 At}}-A)\\
\nonumber&\leq&\frac{CA}{1-2C^2 At},
\end{eqnarray*}
which implies that
$$(z^n)_{n\in \mathbb{N}}\ \
is\ \ uniformly \ \ bounded \ \  in \ \ C([0,T]; B^s_{p,r})\times C([0,T]; B^{s-1}_{p,r}).$$
Using the equations $(T_n)$ and the similar argument in the proof of Lemma 3.1 (1), one can easily prove that
$$\left(\begin{array}{c}
                                \partial_t u^{n+1} \\
                                \partial_t \rho^{n+1} \\
                              \end{array}
                            \right)_{n\in \mathbb{N}}
is\ \ uniformly \ \ bounded \ \  in \ \ C([0,T]; B^{s-1}_{p,r})\times C([0,T]; B^{s-2}_{p,r}).$$
Hence,
$$(z^n)_{n\in \mathbb{N}}\ \
is\ \ uniformly \ \ bounded \ \  in \ \ E^s_{p,r}(T)\times E^{s-1}_{p,r}(T).$$
Now it suffices to show that $(z^n)_{n\in \mathbb{N}}$
is  a Cauchy sequence in $C([0,T]; B^{s-1}_{p,r})\times C([0,T]; B^{s-2}_{p,r})$.
Indeed, For all $m,n\in \mathbb{N}$, from $(T_n)$, we have
\begin{eqnarray*}
&&\partial_t (u^{n+m+1}-u^{n+1})+u^{n+m}\,\partial_x (u^{n+m+1}-u^{n+1})\\
\nonumber&=&P(D)\big(\frac {3} {2} (u^{n+m}-u^n)(u^{n+m}+u^n)+\frac{c}{2}(\rho^{n+m}-\rho^n)(\rho^{n+m}+\rho^n)\big)\\
\nonumber&&+(u^n-u^{n+m})\,\partial_x u^{n+1}
\end{eqnarray*}
and
\begin{eqnarray*}
&&\partial_t (\rho^{n+m+1}-\rho^{n+1})+u^{n+m}\,\partial_x (\rho^{n+m+1}-\rho^{n+1})\\
\nonumber&=&-2\big((\rho^{n+m}-\rho^n)\,\partial_x u^n+\rho^{n+m}\,\partial_x (u^{n+m}-u^n)\big)\\
\nonumber&&+(u^n-u^{n+m})\,\partial_x \rho^{n+1}.
\end{eqnarray*}
Similar to the proof of Lemma 3.1 (1), for $s>\max (2-\frac{1}{p},1+\frac{1}{p},\frac 3 2)$ and $s\neq 2+\frac{1}{p}, 3+\frac{1}{p}$, we can obtain that
$$b^m_{n+1}(t)\leq e^{C U^{n+m}(t)}\big(b^m_{n+1}(0)+C \int_0^t e^{-C U^{n+m}(\tau)}b^m_n(\tau)d^m_n(\tau)d\tau\big),$$
where $b^m_n(t)\triangleq ||(u^{n+m}-u^n)(t)||_{B^{s-1}_{p,r}}+||(\rho^{n+m}-\rho^n)(t)||_{B^{s-2}_{p,r}}$, $U^{n+m}(t)\triangleq \int_0^t ||u^{n+m}(\tau)||_{B^s_{p,r}}d\tau$, and\,\ $d^m_n(t)\triangleq ||u^n(t)||_{B^s_{p,r}}+||u^{n+1}(t)||_{B^s_{p,r}}+||u^{n+m}(t)||_{B^s_{p,r}}+||\rho^n(t)||_{B^{s-1}_{p,r}}+
||\rho^{n+1}(t)||_{B^{s-1}_{p,r}}||\rho^{n+m}(t)||_{B^{s-1}_{p,r}}$.\\
Thanks to Remark 2.1, we have
\begin{eqnarray*}
||\sum\limits_{q=n+1}^{n+m} \Delta_q u_0||_{B^{s-1}_{p,r}}
&=&\big(\sum\limits_{k\geq-1}2^{k(s-1)r}||\Delta_k (\sum\limits_{q=n+1}^{n+m} \Delta_q u_0)||^r_{L^p}\big)^{\frac 1 r}\\
&\leq& C \big(\sum\limits_{k=n}^{n+m+1}2^{-kr}2^{ksr}||\Delta_k{u_0}||^r_{L^p}\big)^{\frac 1 r}\\
&\leq& C2^{-n}||u_0||_{B^s_{p,r}}.
\end{eqnarray*}
Similarly,
$$||\sum\limits_{q=n+1}^{n+m} \Delta_q \rho_0||_{B^{s-2}_{p,r}}\leq C2^{-n}||\rho_0||_{B^{s-1}_{p,r}}.$$
Hence, we obtain
$$b^m_{n+1}(0)\leq C2^{-n}(||u_0||_{B^s_{p,r}}+||\rho_0||_{B^{s-1}_{p,r}}).$$
According to the fact that $(z^n)_{n\in \mathbb{N}}$
is uniformly bounded in $E^s_{p,r}(T)\times E^{s-1}_{p,r}(T)$, we can find a positive constant $C_{T}$ independent of $n,m$ such that
$$b^m_{n+1}(t)\leq C_{T}\big(2^{-n}+\int_0^t b^m_n(\tau)d\tau\big),\quad \forall t\in[0,T].$$
Arguing by induction with respect to the index $n$, we can obtain
\begin{eqnarray*}
b^m_{n+1}(t)
&\leq& C_{T}\big(2^{-n}\sum\limits_{k=0}^n \frac{(2T C_T)^k}{k!}+C^{n+1}_T \int_0^t  \frac{(t-\tau)^n}{n!}d\tau\big)\\
&\leq& \big(C_{T}\sum\limits_{k=0}^n \frac{(2T C_T)^k}{k!}\big)2^{-n}+C_T \frac{(T C_T)^{n+1}}{(n+1)!},
\end{eqnarray*}
which implies the desired result.

On the other hand, for the critical points $s= 2+\frac{1}{p}$ or $3+\frac{1}{p}$, we can apply the interpolation method which has been used in the proof of Lemma 3.1 to show that $(z^n)_{n\in \mathbb{N}}$
is also a  Cauchy sequence in $C([0,T]; B^{s-1}_{p,r})\times C([0,T]; B^{s-2}_{p,r})$ for these two critical cases.
Therefore, we have completed the proof of Lemma 3.2.
\end{proof}

Now we are in the position to prove the main theorem of this section.
\begin{theorem3}
Assume that $1\leq p,r\leq \infty$ and $s>\max (2-\frac{1}{p},1+\frac{1}{p},\frac 3 2)$ with $s\neq2+\frac{1}{p}$. Let
$z_{0}\triangleq
                            \left(\begin{array}{c}
                                u_0 \\
                                \rho_0 \\
                              \end{array}
                            \right)\in B^s_{p,r}\times B^{s-1}_{p,r}$ and
$z\triangleq
\left(\begin{array}{c}
                                u \\
                                \rho \\
                              \end{array}
                            \right)$ be the obtained limit in Lemma 3.2.
Then there exists a time $T>0$ such that $z\in E^s_{p,r}(T)\times E^{s-1}_{p,r}(T)$ is the unique solution to System (3.1), and the mapping $z_0\mapsto z:$ is continuous from
                            $ B^s_{p,r}\times B^{s-1}_{p,r}$ into
$$ C([0,T]; B^{s^{'}}_{p,r})\cap C^1{([0,T]; B^{s^{'}-1}_{p,r})}\!\times\!C([0,T]; B^{s{'}-1}_{p,r})\cap C^1{([0,T]; B^{s^{'}-2}_{p,r})}$$
for all $s^{'}<s$ if $r=\infty$ and $s^{'}=s$ otherwise.
\end{theorem3}

\begin{proof}
We first claim that $z\in E^s_{p,r}(T)\times E^{s-1}_{p,r}(T)$ solves System (3.1).\\
In fact, according to Lemma 3.2 (2) and Proposition 2.2 (7), one can get
$$z\in L^{\infty}([0,T]; B^s_{p,r})\times L^{\infty}([0,T]; B^{s-1}_{p,r}).$$
For all $s^{'}<s$, Lemma 3.2 (2) applied again, together with an interpolation argument yields
$$z^n \rightarrow z,\,\ as \,\ n\to \infty,\,\ in\,\ C([0,T]; B^{s^{'}}_{p,r})\times C([0,T]; B^{s^{'}-1}_{p,r}). $$
Taking limit in $(T_n)$, we can see that
$z$ solves System (3.1) in the sense of $C([0,T]; B^{s^{'}-1}_{p,r})\times C([0,T]; B^{s^{'}-2}_{p,r})$ for all $s^{'}<s$.\\
Making use of the equations in System (1.3) twice and the similar proof in the Claim of Lemma 3.1, together with Lemma 2.1 (4) and Lemma 2.2 yields
$z\in E^s_{p,r}(T)\times E^{s-1}_{p,r}(T).$

On the other hand, the continuity with respect to the initial data in
$$ C([0,T]; B^{s^{'}}_{p,r})\cap C^1{([0,T]; B^{s^{'}-1}_{p,r})}\!\times\!
C([0,T]; B^{s{'}-1}_{p,r})\cap C^1{([0,T]; B^{s^{'}-2}_{p,r})} \quad (\forall\, s^{'}<s)$$
can be obtained by Lemma 3.1 and a simple interpolation argument. While the continuity in $C([0,T]; B^{s}_{p,r})\cap C^1{([0,T]; B^{s-1}_{p,r})}\!\times\!C([0,T]; B^{s-1}_{p,r})\cap C^1{([0,T]; B^{s-2}_{p,r})}$ when $r< \infty$ can be proved through the use of a sequence of viscosity approximation solutions
$\left(\begin{array}{c}
                                u_\varepsilon \\
                                \rho_\varepsilon \\
                              \end{array}
                            \right)_{\varepsilon>0}$
for System (1.3) which converges uniformly in $C([0,T]; B^{s}_{p,r})\cap C^1{([0,T]; B^{s-1}_{p,r})}\!\times\!C([0,T]; B^{s-1}_{p,r})\cap C^1{([0,T]; B^{s-2}_{p,r})}$. This completes the proof of Theorem 3.1.
\end{proof}

\begin{remark3}
(1) Note that for every $s\in\mathbb{R}$, $B^{s}_{2,2}=H^s$. Theorem 3.2 holds true in the corresponding Sobolev spaces with $\frac{3} {2}<s\neq\frac{5} {2}$, which almost improves the result of Theorem 3.1 proved by Kato's theory, where $s\geq 2$ is required. Therefore, Theorem 3.2 together with Theorem 3.1 implies that the conclusion of Theorem 3.1  holds true for all $s>\frac{3} {2}$.\\
(2) As we know, $u_c (t,x)=ce^{-|x-ct|}$ with $c\in \mathbb{R}$ is the solitary wave solution to DP equation \cite{D-H-H}, then the index $s=\frac {3}{2}$ is critical in Besov spaces $B^{s}_{2,r}$ in the following sense \cite{D2}:
System (1.3) is not local well-posedness in $B^{\frac {3}{2}}_{2,\infty}$. More precisely, there exists a global solution $u_1\in L^\infty({\mathbb{R}}^{+}; B^{\frac {3}{2}}_{2,\infty})$ and $v\equiv 0$ to System (1.3) such that for any $T>0$ and $\varepsilon>0$, there exists a solution $u_2\in L^\infty(0,T; B^{\frac {3}{2}}_{2,\infty})$ and $v\equiv 0$ to System (1.3) with
$$||u_1(0)-u_2(0)||_{B^{\frac {3}{2}}_{2,\infty}}\leq \varepsilon\quad but\quad ||u_1-u_2||_{L^\infty(0,T; B^{\frac {3}{2}}_{2,\infty})}\geq 1.$$
\end{remark3}

\section{The precise blow-up scenario}
\newtheorem {remark4}{Remark}[section]
\newtheorem{theorem4}{Theorem}[section]
\newtheorem{lemma4}{Lemma}[section]
\newtheorem{corollary4}{Corollary}[section]

\par
In this section, we will derive the precise blow-up scenario of strong solutions to System (1.3).

Firstly, let us consider the following differential equation:
\begin{equation}
\left\{\begin{array}{ll}q_{t}=u(t,q),\ \ \ \ t\in[0,T), \\
q(0,x)=x,\ \ \ \ x\in\mathbb{R}, \end{array}\right.
\end{equation}
where $u$ denotes the first component of the solution $z$ to
System (1.3).

The following lemmas are very crucial to study the blow-up phenomena of strong solutions to System (1.3).

\begin{lemma4}
\cite{Cf} Let $u\in C([0,T);H^s(\mathbb{R}))\cap
C^1([0,T);H^{s-1}(\mathbb{R})), s\geq 2$. Then Eq.(4.1) has a unique
solution $q\in C^1([0,T)\times \mathbb{R};\mathbb{R})$. Moreover,
the map $q(t,\cdot)$ is an increasing diffeomorphism of $\mathbb{R}$
with
$$
q_{x}(t,x)=\exp\left(\int_{0}^{t}u_{x}(s,q(s,x))ds\right)>0, \ \
\forall(t,x)\in [0,T)\times \mathbb{R}.$$
\end{lemma4}

\begin{lemma4}
Let $z_{0}\triangleq\left(\begin{array}{c}
                                      u_0 \\
                                      \rho_0 \\
                                    \end{array}
                                  \right)\in H^s(\mathbb{R})\times H^{s-1}(\mathbb{R})$ with $s> \frac 3 2$ and $T>0$ be
the maximal existence time of the corresponding solution
$z\triangleq\left(\begin{array}{c}
                                      u \\
                                     \rho\\
                                    \end{array}
                                  \right)$
to System (1.3), which is guaranteed by Remark 3.1 (1). Then we have
\begin{equation}
\rho(t,q(t,x))q^2_{x}(t,x)=\rho_{0}(x),\ \ \ \forall(t,x)\in
[0,T)\times \mathbb{R}.
\end{equation}
Moreover, if there exists a $M>0$ such that $u_{x}(t,x)\geq -M$ for all $(t,x)\in
[0,T)\times \mathbb{R}$, then
$$\|\rho(t,\cdot)\|_{L^{\infty}},\|\rho(t,\cdot)\|_{L^2}\leq e^{2Mt}\|\rho_{0}\|_{H^{s-1}},\ \ \
\forall t\in [0,T).$$
\end{lemma4}

\begin{proof} Differentiating the left-hand side of Eq.(4.2) with
respect to $t$ and making use of (4.1) and System (1.3), we obtain
\begin{eqnarray*}
&&\frac{d}{dt}{(\rho(t,q(t,x))q^2_{x}(t,x))}\\
&=&(\rho_t(t,q)+\rho_x(t,q)q_t(t,x))q^2_{x}(t,x)
\\&&+2\rho(t,q) q_x(t,x) q_{xt}(t,x)\\
&=&\big(\rho_t(t,q)+\rho_x(t,q)u(t,q)+2\rho(t,q)u_x(t,q)\big)q^2_{x}(t,x)\\
&=& 0.
\end{eqnarray*}%
This proves (4.2). By Lemma 4.1, in view of (4.2) and the assumption of the lemma, we obtain for all $t\in[0,T)$
\begin{eqnarray*}
\|\rho(t,\cdot)\|_{L^{\infty}}
&=&\|\rho(t,q(t,\cdot))\|_{L^{\infty}}\\
&=&\|e^{-2\int_0^t u_x(s,\cdot)ds}\rho_0(\cdot)\|_{L^{\infty}}\\
&\leq&e^{2Mt}\|\rho_0(\cdot)\|_{L^{\infty}}.
\end{eqnarray*}
By (4.2) and Lemma 4.1, we get
\begin{eqnarray*}
\int_{\mathbb{R}}{|\rho(t,x)|}^2 dx
&=&\int_{\mathbb{R}}{|\rho(t,q(t,x))|}^2 q_x(t,x)dx =\int_{\mathbb{R}}{|\rho_0(x)|}^2 q_x^{-3}(t,x)dx\\
&\leq&e^{3Mt}\int_{R}{|\rho_{0}(x)|}^2 dx, \ \ \ \forall t\in [0,T).
\end{eqnarray*}
 This completes the proof of the lemma.
\end{proof}

As mentioned in the Introduction, the $H^1$-norm of the solutions to DP equation is not conserved. However, what saves the game in some sense is to establish a priori estimate for the $L^\infty$-norm of the first component $u$ of the strong solutions to System (1.3).
\begin{lemma4}
Let $z_{0}=\left(
                                                     \begin{array}{c}
                                                       u_{0} \\
                                                       \rho_{0} \\
                                                     \end{array}
                                                   \right)
\in H^s(\mathbb{R})\times H^{s-1}(\mathbb{R})$ with $s>\frac 3 2$ and
$T$ be the maximal existence time of the solution $z=\left(
                                                     \begin{array}{c}
                                                       u \\
                                                      \rho \\
                                                     \end{array}
                                                   \right) $ to System (1.3), which is guaranteed by Remark 3.1 (1).
Assume that there is a $M>0$ such that
$\|\rho(t,\cdot)\|_{L^{\infty}},\|\rho(t,\cdot)\|_{L^2}\leq e^{2Mt}\|\rho_{0}\|_{H^{s-1}}$ for all $t\in [0,T)$.
Then for all $t\in [0,T)$, we have
\begin{eqnarray}
\|u(t)\|^2_{L^2}\leq
2\,e^{2|c|t}\big(2\|u_0\|^2_{L^2}+|c|\,t(1+8Mt)(e^{2Mt}\|\rho_{0}\|_{H^{s-1}})^4\big)
\end{eqnarray}
and
\begin{eqnarray}
&&\ \ \ \ ||u(t)||_{L^{\infty}} \\
&\leq&\frac{3}{2}\,t\,e^{2|c|t}\big(2\|u_0\|^2_{L^2}+|c|\,t(1+8Mt)(e^{2Mt}\|\rho_{0}\|_{H^{s-1}})^{4}
+\frac{|c|}{4}(e^{2Mt}\|\rho_{0}\|_{H^{s-1}})^{2}\big) \nonumber \\
&&+\|u_0\|_{L^{\infty}}\nonumber \\
&\triangleq& J(t).\nonumber
\end{eqnarray}
\end{lemma4}

\begin{proof}
By a standard density argument, here we may assume $s\geq 3$ to prove the lemma.
Set $w\triangleq(4-\partial_x^2)^{-1}u$. By the first equation of System (1.1) and the fact that $(\hat{m_t},\hat w)=(\hat m,\hat{w_t})$ or $\int_{\mathbb{R}}m_twdx=\int_{\mathbb{R}}m w_tdx$, we have
\begin{eqnarray*}
\frac{1}{2}\frac{d}{dt}\int_{\mathbb{R}}m w dx
&=&\frac{1}{2}\int_{\mathbb{R}}m_t w dx+
\frac{1}{2}\int_{\mathbb{R}}m w_tdx=\int_{\mathbb{R}}m_twdx\\
&=& -3\int_{\mathbb{R}}w m u_x dx-\int_{\mathbb{R}}w m_x u dx-c\int_{\mathbb{R}}w\rho\rho_x dx\\
&=& -\int_{\mathbb{R}}w(mu)_x dx-2\int_{\mathbb{R}}w m u_x dx+\frac{c}{2}\int_{\mathbb{R}}w_x\rho^2dx.
\end{eqnarray*}
While
\begin{eqnarray*}
\int_{\mathbb{R}}w(mu)_x dx=-\int_{\mathbb{R}}w_x mu dx =\int_{\mathbb{R}}w_xu^2dx-\int_{\mathbb{R}}w_x u_x^2 dx,
\end{eqnarray*}
and
\begin{eqnarray*}
2\int_{\mathbb{R}}w mu_x dx=-\int_{\mathbb{R}}w_xu^2
+\int_{\mathbb{R}}w_x u_x^2 dx.
\end{eqnarray*}
Combining the above three equalities, we deduce that
\begin{eqnarray*}
\frac{d}{dt}\int_{\mathbb{R}}m w dx=c\int_{\mathbb{R}}w_x\rho^2dx.
\end{eqnarray*}
Integrating from $0$ to $t$ on both sides of the above equality, we have
\begin{eqnarray*}
\int_{\mathbb{R}}m w dx=\int_{\mathbb{R}}m_0w_0dx+c\int_0^t\int_{\mathbb{R}}w_x\rho^2dxds,
\end{eqnarray*}
which implies
\begin{eqnarray*}
\|u(t)\|^2_{L^2}
&=&\|\hat{u}(t)\|^2_{L^2}\leq4\int_{\mathbb{R}}\frac{1+\xi^2}{4+\xi^2}|\hat{u}(t,\xi)|^2d\xi=4(\hat{m}(t),\hat{w}(t))\\
&=&4(m(t),w(t))=4(m_0,w_0)+4c\int_0^t\int_{\mathbb{R}}w_x\rho^2dxds\\
&\leq&4\|u_0\|^2_{L^2}+4c\int_0^t\int_{\mathbb{R}}w_x\rho^2dxds.
\end{eqnarray*}
Note that
\begin{eqnarray*}
\|w_x(t)\|^2_{L^2}
&=&\|\partial_x(4-\partial^2_x)^{-1}u(t)\|^2_{L^2}\\
&\leq& \|u(t)\|^2_{L^2}.
\end{eqnarray*}
Besides, by the assumption of the lemma,  we have
\begin{eqnarray*}
\|\rho(t,\cdot)\|^4_{L^4}\leq \|\rho(t,\cdot)\|^2_{L^{\infty}}\|\rho(t,\cdot)\|^2_{L^2}\leq (e^{2Mt}\|\rho_{0}\|_{H^{s-1}})^4,\quad
\forall\,\ t\in[0,T).
\end{eqnarray*}
Hence,
\begin{eqnarray*}
\|u(t)\|^2_{L^2}
&\leq& 4\|u_0\|^2_{L^2}+2c\int_0^t(\|w_x(s,\cdot)\|^2_{L^2}+\|\rho(s,\cdot)\|^4_{L^4})d s\\
&\leq& 4\|u_0\|^2_{L^2}+2|c|t(e^{2Mt}\|\rho_{0}\|_{H^{s-1}})^4+ 2|c|\int_0^t\|u(s)\|^2_{L^2}d s.
\end{eqnarray*}
By Gronwall's inequality, we can reach (4.3).

Next we prove (4.4). Indeed, by the first equation in System (1.3), we have
$$u_t+u u_x= -\partial_x p\ast (\frac 3 2 {u^2}+\frac {c}{2} {\rho^2}).$$
Applying Young's inequality and noting that $\|\partial_x p \|_{L^\infty}\leq\frac{1}{2}$, we have
\begin{eqnarray*}
&&\|-\partial_xp\ast
(\frac{3}{2}u^2+\frac{c}{2}\rho^2)\|_{L^\infty}\\&\leq&
\|\partial_xp\|_{L^\infty}\|\frac{3}{2}u^2+\frac{c}{2}\rho^2\|_{L^1}\\
&\leq&\frac{3}{4}\|u\|^2_{L^2}+\frac{|c|}{4}\|\rho\|^2_{L^2}\\
&\leq&\frac{3}{4}\|u\|^2_{L^2}+\frac{|c|}{4}(e^{2Mt}\|\rho_{0}\|_{H^{s-1}})^2.
\end{eqnarray*}
Besides, in view of (4.1), we have
\begin{eqnarray*}
\frac{du(t,q(t,x))}{dt}=u_t(t,q(t,x))+u_x(t,q(t,x))q_t(t,x)=(u_t+uu_x)(t,q(t,x)).
\end{eqnarray*}
Thanks to (4.3) and the facts above, we deduce
\begin{eqnarray*}
-P(t)\leq\frac{du(t,q(t,x))}{dt}\leq P(t),
\end{eqnarray*}
where
\begin{eqnarray*}
P(t)\triangleq\frac{3}{2}e^{2|c|t}\big(2\|u_0\|^2_{L^2}+|c|\,t(1+8Mt)(e^{2Mt}\|\rho_{0}\|_{H^{s-1}})^{4}+\frac{|c|}{4}(e^{2Mt}\|\rho_{0}\|_{H^{s-1}})^{2}\big).
\end{eqnarray*}
Integrating the above inequalities with respect to $t<T$ on $[0,t]$ yields
\begin{eqnarray*}
-t\,P(t)+u_0(x)\leq u(t,q(t,x)) \leq t\,P(t)+u_0(x).
\end{eqnarray*}
Therefore, in view of Lemma 4.1, we get the desired result. This completes the proof of the lemma.
\end{proof}

\begin{corollary4}
Let $z_{0}=\left(
                                                     \begin{array}{c}
                                                       u_{0} \\
                                                       \rho_{0} \\
                                                     \end{array}
                                                   \right)
\in H^s(\mathbb{R})\times H^{s-1}(\mathbb{R})$ with $s>\frac 3 2$ and
$T$ be the maximal existence time of the solution $z=\left(
                                                     \begin{array}{c}
                                                       u \\
                                                      \rho \\
                                                     \end{array}
                                                   \right) $ to System (1.3), which is guaranteed by Remark 3.1 (1).
If $\partial_x u \in L^1(0,T; {L^\infty})$,
then for all $t\in [0,T)$, we have
\begin{eqnarray*}
&&\ \ \ \ ||u(t)||_{L^{\infty}} \\
&\leq&\frac{3}{2}\,t\,e^{2|c|t}\big(2\|u_0\|^2_{L^2}+|c|\,t(1+8\int_0^t ||\partial_x u(\tau)||_{L^\infty} d\tau)(e^{2\int_0^t ||\partial_x u(\tau)||_{L^\infty} d\tau}\|\rho_{0}\|_{H^{s-1}})^{4} \nonumber \\
&&+\frac{|c|}{4}(e^{2\int_0^t ||\partial_x u(\tau)||_{L^\infty} d\tau}\|\rho_{0}\|_{H^{s-1}})^{2}\big)+\|u_0\|_{L^{\infty}}\nonumber \\
&\triangleq& L(t).\nonumber
\end{eqnarray*}
\end{corollary4}

\begin{proof}
By the proof of Lemma 4.2, we also have
\begin{eqnarray*}
\|\rho(t,\cdot)\|_{L^{\infty}},\,\ \|\rho(t,\cdot)\|_{L^2}
\leq e^{2\int_0^t ||\partial_x u(\tau)||_{L^\infty} d\tau} \|\rho_0\|_{H^{s-1}},\quad  \forall \,\ t\in [0,T).
\end{eqnarray*}
It is then easy to prove the corollary by a similar argument as in the proof of Lemma 4.3.
\end{proof}

\begin{theorem4}
Let $z_{0}=\left(
                                                     \begin{array}{c}
                                                       u_{0} \\
                                                       \rho_{0} \\
                                                     \end{array}
                                                   \right)
\in H^s(\mathbb{R})\times H^{s-1}(\mathbb{R})$ with $s>\frac 3 2$ and
$T$ be the maximal existence time of the solution $z=\left(
                                                     \begin{array}{c}
                                                       u \\
                                                      \rho \\
                                                     \end{array}
                                                   \right) $ to System (1.3), which is guaranteed by Remark 3.1 (1).
If $T<\infty$, then $$\int_0^T ||\partial_x u(\tau)||_{L^\infty} d\tau=\infty.$$
\end{theorem4}

\begin{proof}
 We will prove the theorem by induction with respect to the regular index $s$ $(s>\frac 3 2)$ as follows.

Step 1. For $s\in(\frac{3}{2},2)$, by Lemma 2.3 and the second equation of System (1.3), we have
\begin{eqnarray*}
||\rho(t)||_{H^{s-1}}&\leq& ||\rho_0||_{H^{s-1}}\,+\, \ C \int_0^t ||\partial_x u(\tau)\, \rho(\tau)||_{H^{s-1}}d\tau\\
 &&+C\int_0^t ||\rho(\tau)||_{H^{s-1}}(||u(\tau)||_{L^{\infty}}+||\partial_x u(\tau)||_{L^{\infty}}) d\tau.
\end{eqnarray*}
Applying (2.1), we get
\begin{eqnarray}
||\partial_x u\, \rho||_{H^{s-1}}\leq
C(||\partial_x u||_{H^{s-1}}||\rho||_{L^\infty}+||\partial_x u||_{L^\infty} ||\rho||_{H^{s-1}}).
\end{eqnarray}
Thus,
\begin{eqnarray}
||\rho(t)||_{H^{s-1}}&\leq& ||\rho_0||_{H^{s-1}}\,+\, \ C \int_0^t ||\partial_x u(\tau)||_{H^{s-1}}||\rho(\tau)||_{L^\infty}d\tau\\
 &&+C\int_0^t ||\rho(\tau)||_{H^{s-1}}(||u(\tau)||_{L^{\infty}}+||\partial_x u(\tau)||_{L^{\infty}}) d\tau.\nonumber
\end{eqnarray}
On the other hand, thanks to Lemma 2.1 (3) and the first equation of System (1.3), we have ($\forall s>1$, indeed)
\begin{eqnarray*}
||u(t)||_{H^s}&\leq& ||u_0||_{H^s}\,+\, \ C \int_0^t ||P(D)(\frac 3 2 {u^2}+\frac {c}{2} {\rho^2})(\tau)||_{H^s}d\tau\\
 &&+C\int_0^t ||u(\tau)||_{H^s}||\partial_x u(\tau)||_{L^{\infty}} d\tau.
\end{eqnarray*}
By Proposition 2.2 (8) and (2.1), we have
\begin{eqnarray*}
&& ||P(D)(\frac 3 2 {u^2}+\frac {c}{2} {\rho^2})||_{H^s}\\
&\leq& C||\frac 3 2 {u^2}+\frac {c}{2} {\rho^2}||_{H^{s-1}}\nonumber \\
&\leq& C(||u||_{H^{s-1}}||u||_{L^\infty}+||\rho||_{H^{s-1}}||\rho||_{L^\infty}).\nonumber
\end{eqnarray*}
Hence,
\begin{eqnarray}
\quad\quad ||u(t)||_{H^s}&\leq& ||u_0||_{H^s}\,+\, \ C \int_0^t ||u(\tau)||_{H^s}(||u(\tau)||_{L^{\infty}}+||\partial_x u(\tau)||_{L^{\infty}})d\tau\\
&&+C\int_0^t ||\rho(\tau)||_{H^{s-1}}||\rho(\tau)||_{L^{\infty}} d\tau. \nonumber
\end{eqnarray}
Combining (4.6) and (4.7), we obtain
\begin{eqnarray*}
||u(t)||_{H^s}+||\rho(t)||_{H^{s-1}}&\leq& ||u_0||_{H^s}\,+\,||\rho_0||_{H^{s-1}}+C \int_0^t (||u||_{H^s}+||\rho||_{H^{s-1}})\\
&&\times (||u||_{L^{\infty}}+||\partial_x u||_{L^{\infty}}+||\rho||_{L^{\infty}})d\tau.
\end{eqnarray*}
Thanks to Gronwall's inequality, we have
\begin{eqnarray}
\quad\quad ||u(t)||_{H^s}+||\rho(t)||_{H^{s-1}}&\leq& (||u_0||_{H^s}\,+\,||\rho_0||_{H^{s-1}})\\
&&\times e^{C \int_0^t (||u(\tau)||_{L^{\infty}}+||\partial_x u(\tau)||_{L^{\infty}}+||\rho(\tau)||_{L^{\infty}})d\tau}.\nonumber
\end{eqnarray}
Therefore, if $T<\infty$ satisfies $\int_0^T ||\partial_x u(\tau)||_{L^\infty} d\tau<\infty$, then we deduce from (4.8), Corollary 4.1 and the fact $\|\rho(t,\cdot)\|_{L^{\infty}}\leq e^{2\int_0^t ||\partial_x u(\tau)||_{L^\infty} d\tau} \|\rho_0\|_{H^{s-1}}$ that
\begin{eqnarray}
&& ||u(t)||_{H^s}+||\rho(t)||_{H^{s-1}}\\
&\leq& (||u_0||_{H^s}\,+\,||\rho_0||_{H^{s-1}})\nonumber \\
&&\times e^{Ct(e^{2{\int_0^t ||\partial_x u(\tau)||_{L^\infty} d\tau}}\|\rho_{0}\|_{H^{s-1}}+L(t))+C \int_0^t ||\partial_x u(\tau)||_{L^{\infty}}d\tau}.\nonumber
\end{eqnarray}
Hence,
\begin{eqnarray}
\limsup\limits_{t\to T}(||u(t)||_{H^s}+||\rho(t)||_{H^{s-1}})<\infty,
\end{eqnarray}
which contradicts the assumption that $T<\infty$ is the maximal existence time. This completes the proof of the theorem for $s\in(\frac 3 2, 2)$.

Step 2. For $s\in[2,\frac 5 2)$, Lemma 2.1 (1) applied to the second equation of System (1.3), we get
\begin{eqnarray*}
||\rho(t)||_{H^{s-1}}&\leq& ||\rho_0||_{H^{s-1}}\,+\, \ C \int_0^t ||\partial_x u(\tau)\, \rho(\tau)||_{H^{s-1}}d\tau\\
 &&+C\int_0^t ||\rho(\tau)||_{H^{s-1}}||\partial_x u(\tau)||_{H^{\frac 1 2}\cap L^\infty} d\tau,
\end{eqnarray*}
together with (4.5) implies that
\begin{eqnarray}
||\rho(t)||_{H^{s-1}}&\leq& ||\rho_0||_{H^{s-1}}\,+\, \ C \int_0^t ||\partial_x u(\tau)||_{H^{s-1}}||\rho(\tau)||_{L^\infty}d\tau\\
 &&+C\int_0^t ||\rho(\tau)||_{H^{s-1}}||\partial_x u(\tau)||_{H^{\frac 1 2}\cap L^\infty} d\tau.\nonumber
\end{eqnarray}
By (4.7) and (4.11), we have
\begin{eqnarray*}
||u(t)||_{H^s}+||\rho(t)||_{H^{s-1}}&\leq& ||u_0||_{H^s}\,+\,||\rho_0||_{H^{s-1}}+C \int_0^t (||u||_{H^s}+||\rho||_{H^{s-1}})\\
&&\times (||u||_{H^{\frac 3 2 +\varepsilon_0}}+||\rho||_{L^{\infty}})d\tau,
\end{eqnarray*}
where $\varepsilon_0\in (0,\frac 1 2)$ and we used the fact that $H^{\frac 1 2 +\varepsilon_0}\hookrightarrow H^{\frac 1 2}\cap L^\infty$.
Thanks to Gronwall's inequality again, we have
\begin{eqnarray}
\quad\quad ||u(t)||_{H^s}+||\rho(t)||_{H^{s-1}}&\leq& (||u_0||_{H^s}\,+\,||\rho_0||_{H^{s-1}})\\
&&\times e^{C \int_0^t (||u(\tau)||_{H^{\frac 3 2 +\varepsilon_0}}+||\rho(\tau)||_{L^{\infty}})d\tau}.\nonumber
\end{eqnarray}
Therefore if $T<\infty$ satisfies $\int_0^T ||\partial_x u(\tau)||_{L^\infty} d\tau<\infty$, then we deduce from (4.12), (4.10)  with $\frac 3 2 +\varepsilon_0\in(\frac 3 2,2)$ and the fact $\|\rho(t,\cdot)\|_{L^{\infty}}\leq e^{2\int_0^t ||\partial_x u(\tau)||_{L^\infty} d\tau} \|\rho_0\|_{H^{s-1}}$ that
\begin{eqnarray}
\limsup\limits_{t\to T}(||u(t)||_{H^s}+||\rho(t)||_{H^{s-1}})<\infty,
\end{eqnarray}
which contradicts the assumption that $T<\infty$ is the maximal existence time. This completes the proof of the theorem for $s\in[2,\frac 5 2)$.

Step 3. For $s\in(2,3)$, by differentiating the second equation of System (1.3) with respect to $x$, we have
\begin{eqnarray*}
\partial_t \rho_x+ u\,\partial_x \rho_x+ 3u_x \rho_x+ 2u_{xx}\rho=0.
\end{eqnarray*}
By Lemma 2.3, we get
\begin{eqnarray*}
||\partial_x \rho(t)||_{H^{s-2}}&\leq& ||\partial_x \rho_0||_{H^{s-2}}\,+\, \ C \int_0^t ||(3u_x \rho_x+ 2\rho u_{xx})(\tau)||_{H^{s-2}}d\tau\\
&&+C\int_0^t ||\partial_x \rho(\tau)||_{H^{s-2}}(||u(\tau)||_{L^{\infty}}+||\partial_x u(\tau)||_{L^{\infty}}) d\tau.\nonumber
\end{eqnarray*}
Thanks to (2.3), we have
\begin{eqnarray*}
||u_x \rho_x||_{H^{s-2}}\leq
C(||\partial_x u||_{H^{s-1}}||\rho||_{L^\infty}+||\partial_x u||_{L^\infty} ||\partial_x \rho||_{H^{s-2}})
\end{eqnarray*}
and
\begin{eqnarray*}
||\rho u_{xx}||_{H^{s-2}}\leq
C(||\rho||_{H^{s-1}}||\partial_x u||_{L^\infty}+||\rho||_{L^\infty} ||u_{xx}||_{H^{s-2}}).
\end{eqnarray*}
Hence,
\begin{eqnarray*}
||\partial_x \rho(t)||_{H^{s-2}}&\leq& ||\partial_x \rho_0||_{H^{s-2}}\,+\, \ C \int_0^t (||u(\tau)||_{H^s}+||\rho(\tau)||_{H^{s-1}})\\
&&\times (||u(\tau)||_{L^{\infty}}+||\partial_x u(\tau)||_{L^{\infty}}+||\rho(\tau)||_{L^{\infty}}) d\tau,\nonumber
\end{eqnarray*}
which together with (4.7) and (4.6) with $s-2$ instead of $s-1$, yields that
\begin{eqnarray*}
||u(t)||_{H^s}+||\rho(t)||_{H^{s-1}}&\leq& ||u_0||_{H^s}\,+\,||\rho_0||_{H^{s-1}}+C \int_0^t (||u||_{H^s}+||\rho||_{H^{s-1}})\\
&&\times (||u||_{L^{\infty}}+||\partial_x u||_{L^{\infty}}+||\rho||_{L^{\infty}})d\tau.
\end{eqnarray*}
Similar to Step 1, we can easily prove the theorem for $s\in(2,3)$.

Step 4. For $s=k\in\mathbb{N}$ and $k\geq3$, by differentiating the second equation of System (1.3) $k-2$ times with respect to $x$, we get
\begin{eqnarray*}
(\partial_t+u \partial_x) \partial^{k-2}_x \rho + \sum\limits_{l_1+l_2=k-3,l_1,l_2\geq 0} {C_{l_1,l_2}{\partial^{l_1+1}_x u}\,{\partial^{l_2+1}_x \rho}}+2\rho\, \partial^{k-1}_x u=0,
\end{eqnarray*}
which together with Lemma 2.1 (1), implies that
\begin{eqnarray*}
||\partial^{k-2}_x \rho(t)||_{H^1}&\leq& ||\partial^{k-2}_x \rho_0||_{H^1}+C \int_0^t ||\partial^{k-2}_x \rho(\tau)||_{H^1}||\partial_x u(\tau)||_{H^{\frac 1 2}\cap L^\infty} d\tau\\
&+&C \int_0^t ||\sum\limits_{l_1+l_2=k-3,l_1,l_2\geq 0} {C_{l_1,l_2}{\partial^{l_1+1}_x u}\,{\partial^{l_2+1}_x \rho}}+2\rho\, \partial^{k-1}_x u||_{H^1} d\tau.
\end{eqnarray*}
Since $H^1$ is an algebra, it follows that
\begin{eqnarray*}
||\rho\, \partial^{k-1}_x u||_{H^1}\leq C ||\rho||_{H^1}||\partial^{k-1}_x u||_{H^1}\leq C||\rho||_{H^1}||u||_{H^s}
\end{eqnarray*}
and
\begin{eqnarray*}
||\sum\limits_{l_1+l_2=k-3,l_1,l_2\geq 0} {C_{l_1,l_2}{\partial^{l_1+1}_x u}\,{\partial^{l_2+1}_x \rho}}||_{H^1}\leq C ||u||_{H^{s-1}}||\rho||_{H^{s-1}}.
\end{eqnarray*}
Then, we have
\begin{eqnarray}
||\partial^{k-2}_x \rho(t)||_{H^1}&\leq& ||\partial^{k-2}_x \rho_0||_{H^1}+C \int_0^t(||u(\tau)||_{H^s}+||\rho(\tau)||_{H^{s-1}})\\
&&\times (||u(\tau)||_{H^{s-1}}+||\rho(\tau)||_{H^1})d\tau.\nonumber
\end{eqnarray}
By the classical Gagliardo-Nirenberg  inequality, we have for $\sigma\in(0,1)$,
\begin{eqnarray*}
||\rho(t)||_{H^{s-1}}\leq C (||\rho(t)||_{H^{\sigma}}+||\partial^{k-2}_x \rho(t)||_{H^1}),
\end{eqnarray*}
which together with (4.14), (4.7) and (4.6) with $\sigma$ instead of $s-1$, yields that
\begin{eqnarray*}
||u(t)||_{H^s}+||\rho(t)||_{H^{s-1}}&\leq& C(||u_0||_{H^s}\,+\,||\rho_0||_{H^{s-1}})+C \int_0^t (||u||_{H^s}+||\rho||_{H^{s-1}})\\
&&\times (||u||_{H^{s-1}}+||\rho||_{H^1})d\tau.\nonumber
\end{eqnarray*}
By Gronwall's inequality, we obtain
\begin{eqnarray}
||u(t)||_{H^s}+||\rho(t)||_{H^{s-1}}&\leq& C(||u_0||_{H^s}\,+\,||\rho_0||_{H^{s-1}})\\
&&\times e^{C\int_0^t (||u(\tau)||_{H^{s-1}}+||\rho(\tau)||_{H^1})d\tau}.\nonumber
\end{eqnarray}
If $T<\infty$ satisfies $\int_0^T ||\partial_x u(\tau)||_{L^\infty} d\tau<\infty$, applying Step 3 and arguing by induction assumption, we can obtain that $||u(t)||_{H^{s-1}}+||\rho(t)||_{H^1}$ is uniformly bounded. Thanks to (4.15), we get
\begin{eqnarray}
\limsup\limits_{t\to T}(||u(t)||_{H^s}+||\rho(t)||_{H^{s-1}})<\infty,
\end{eqnarray}
which contradicts the assumption that $T<\infty$ is the maximal existence time. This completes the proof of the theorem for $s=k\in\mathbb{N}$ and $k\geq3$.

Step 5. For $s\in(k,k+1),\,k\in\mathbb{N} $ and $k\geq 3$, by differentiating the second equation of System (1.3) $k-1$ times with respect to $x$,
we get
\begin{eqnarray*}
(\partial_t+u \partial_x) \partial^{k-1}_x \rho + \sum\limits_{l_1+l_2=k-2,l_1,l_2\geq 0} {C_{l_1,l_2}{\partial^{l_1+1}_x u}\,{\partial^{l_2+1}_x \rho}}+2\rho\, \partial^{k}_x u=0.
\end{eqnarray*}
Applying Lemma 2.3 with $s-k\in(0,1)$, we have
\begin{eqnarray}
&&\ \ \ \||\partial^{k-1}_x \rho(t)||_{H^{s-k}}\\
&\leq& ||\partial^{k-1}_x \rho_0||_{H^{s-k}}+C \int_0^t ||\partial^{k-1}_x \rho(\tau)||_{H^{s-k}}
(||u(\tau)||_{L^\infty}+||\partial_x u(\tau)||_{L^\infty}) d\tau \nonumber\\
&&+C \int_0^t ||\sum\limits_{l_1+l_2=k-2,l_1,l_2\geq 0} {C_{l_1,l_2}{\partial^{l_1+1}_x u(\tau)}\,{\partial^{l_2+1}_x \rho(\tau)}}+2\rho(\tau)\, \partial^{k}_x u(\tau)||_{H^1} d\tau.\nonumber
\end{eqnarray}
For each $\varepsilon_0\in(0,\frac 1 2)$, using (2.3) and the fact that $H^{{\frac {1} {2}}+\varepsilon_0}\hookrightarrow L^\infty$, we have
\begin{eqnarray}
||\rho \partial^{k}_x u||_{H^{s-k}}
&\leq& C(||\partial^{k}_x u||_{H^{s-k}}||\rho||_{L^\infty}+||\partial^{k-1}_x u||_{L^\infty} ||\rho||_{H^{s-k+1}}) \\
&\leq& C(||u||_{H^s}||\rho||_{L^\infty}+||u||_{H^{k-\frac {1} {2}+\varepsilon_0}} ||\rho||_{H^{s-k+1}})\nonumber
\end{eqnarray}
and
\begin{eqnarray}
&&\ \ \ \ ||\sum\limits_{l_1+l_2=k-2,l_1,l_2\geq 0} {C_{l_1,l_2}{\partial^{l_1+1}_x u}\,{\partial^{l_2+1}_x \rho}}||_{H^{s-k}}\\
&\leq& C\sum\limits_{l_1+l_2=k-2,l_1,l_2\geq 0} {C_{l_1,l_2}(||\partial^{l_1+1}_x u||_{L^\infty}||\partial^{l_2+1}_x \rho}||_{H^{s-k}}\nonumber\\
&&+||\partial^{l_1+1}_x u||_{H^{s-k+1}}||\partial^{l_2}_x \rho||_{L^\infty})\nonumber\\
&\leq& C(||u||_{H^s}||\rho||_{H^{k-\frac {3} {2}+\varepsilon_0}}+||u||_{H^{k-\frac {1} {2}+\varepsilon_0}} ||\rho||_{H^{s-1}}).\nonumber
\end{eqnarray}
Combining (4.17), (4.18) and (4.19), we can get
\begin{eqnarray*}
||\partial^{k-1}_x \rho(t)||_{H^{s-k}}
&\leq& ||\partial^{k-1}_x \rho_0||_{H^{s-k}}+C \int_0^t (||u(\tau)||_{H^s}+||\rho(\tau)||_{H^{s-1}})\\
&&\times (||u||_{H^{k-\frac {1} {2}+\varepsilon_0}}+||\rho||_{H^{k-\frac {3} {2}+\varepsilon_0}})d\tau,
\end{eqnarray*}
which together with (4.7) and (4.6) with $s-k\in(0,1)$ instead of $s-1$, yields that
\begin{eqnarray*}
||u(t)||_{H^s}+||\rho(t)||_{H^{s-1}}&\leq& C(||u_0||_{H^s}\,+\,||\rho_0||_{H^{s-1}})+C \int_0^t (||u||_{H^s}+||\rho||_{H^{s-1}})\\
&&\times (||u||_{H^{k-\frac {1} {2}+\varepsilon_0}}+||\rho||_{H^{k-\frac {3} {2}+\varepsilon_0}})d\tau.
\end{eqnarray*}
Thanks to Gronwall's inequality again, we obtain
\begin{eqnarray*}
||u(t)||_{H^s}+||\rho(t)||_{H^{s-1}}&\leq& C(||u_0||_{H^s}\,+\,||\rho_0||_{H^{s-1}})\\
&&\times e^{C\int_0^t (||u||_{H^{k-\frac {1} {2}+\varepsilon_0}}+||\rho||_{H^{k-\frac {3} {2}+\varepsilon_0}})d\tau}.\nonumber
\end{eqnarray*}
Noting that $k-\frac {1} {2}+\varepsilon_0 <k$, $k-\frac {3} {2}+\varepsilon_0 <k-1$ and $k\geq 3$, and applying Step 3 and the similar argument by induction as in Step 4, we can easily get the desired result.

Consequently, we have completed the proof of the theorem from Step 1 to Step 5.
\end{proof}

The following main theorem of this section shows the precise blow-up scenario for
sufficiently regular solutions to System (1.3).

\begin{theorem4} Let $z_{0}=\left(\begin{array}{c}
                                      u_0 \\
                                      \rho_0 \\
                                    \end{array}
                                  \right)\in H^s(\mathbb{R})\times H^{s-1}(\mathbb{R})$
with $s> \frac 3 2$ and $T>0$ be the maximal existence time of the corresponding solution
$z=\left(\begin{array}{c}
                                      u \\
                                      \rho \\
                                    \end{array}
                                  \right)$
to System (1.3), which is guaranteed by Remark 3.1 (1). Then the corresponding solution $z$ blows up in finite time
if and only if
$$\liminf_{t\rightarrow T}\inf_{x\in\mathbb{R}}\{u_{x}(t,x)\}=-\infty.$$
\end{theorem4}

\begin{proof}
As mentioned earlier, we only need to prove the theorem for $s\geq 3$. Assume that the solution $z$ blows up in finite time ($T<\infty$) and there exists a $M>0$ such that
\begin{eqnarray}
u_{x}(t,x)\geq -M,\ \ \ \forall(t,x)\in [0,T)\times\mathbb{R}.
\end{eqnarray}
By Lemma 4.2, we have
$$\|\rho(t,\cdot)\|_{L^{\infty}},\|\rho(t,\cdot)\|_{L^2}\leq e^{2Mt}\|\rho_{0}\|_{H^{s-1}},\,\ \forall\,  t\in [0,T).$$
Differentiating the first equation in System (1.3)
with respect to $x$ and noting that $\partial^2_x p\ast f=p\ast f-f$, we have
\begin{eqnarray}
u_{tx}=-u_x^2-u u_{xx}-p\ast(\frac{3}{2}u^2+\frac{c}{2}\rho^2)+\frac{3}{2}u^2+\frac{c}{2}\rho^2.
\end{eqnarray}
Note that
\begin{eqnarray}
\frac{du_x(t,q(t,x))}{dt}
&=&u_{xt}(t,q(t,x))+u_{xx}(t,q(t,x))q_t(t,x)\\
&=&(u_{tx}+u u_{xx})(t,q(t,x))\nonumber.
\end{eqnarray}
By (4.21) and (4.22), in view of $u^2_x\geq 0$, $ p\ast u^2 \geq 0$,
$\|p\ast \rho^2\|_{L^{\infty}}\leq \|p\|_{L^{1}}\|\rho\|^2_{L^{\infty}}\leq (e^{2Mt}\|\rho_{0}\|_{H^{s-1}})^2$ and (4.4),
we obtain
\begin{eqnarray*}
&&\ \ \ \
\frac{du_x(t,q(t,x))}{dt}\\
\nonumber&=&-u^2_{x}(t,q(t,x))-p\ast(\frac{3}{2}u^2+\frac{c}{2}\rho^2)(t,q(t,x))+(\frac{3}{2}u^2+\frac{c}{2}\rho^2)(t,q(t,x))\\
\nonumber&\leq&|c|(e^{2Mt}\|\rho_{0}\|_{H^{s-1}})^{2}+\frac{3}{2}J^2(t).
\end{eqnarray*}
Integrating the above inequality with respect to $t<T$ on $[0,t]$ yields that
$$u_x(t,q(t,x))\leq  u_x(0)+|c|t(e^{2Mt}\|\rho_{0}\|_{H^{s-1}})^{2}+\frac{3}{2}t J^2(t),\quad \forall\,t\in[0,T).$$
Then for all $t\in[0,T)$, we have
\begin{eqnarray*}
\sup\limits_{x\in\mathbb{R}} u_x (t,x)
&\leq& ||\partial_x u_0||_{L^\infty}+|c|t(e^{2Mt}\|\rho_{0}\|_{H^{s-1}})^{2}+\frac{3}{2}t J^2(t)\\
&\leq& ||u_{0}||_{H^s}+|c|t(e^{2Mt}\|\rho_{0}\|_{H^{s-1}})^{2}+\frac{3}{2}t J^2(t),
\end{eqnarray*}
which together with (4.20) and $T<\infty$, implies that
$$\int_0^T ||\partial_x u(\tau)||_{L^\infty} d\tau<\infty.$$
This contradicts Theorem 4.1.

On the other hand, by Sobolev's imbedding theorem, we can see that if
$$\liminf_{t\rightarrow T}\inf_{x\in\mathbb{R}}\{u_{x}(t,x)\}=-\infty ,$$
then the solution $z$ will blow up in finite time. This completes the proof of the theorem.
\end{proof}

\begin{remark4}
Theorem 4.2 implies that the blow-up phenomena of the solution $z$ to System (1.3) depends only on the slope of the first component $u$. That is, the first component $u$ must blow up before the second component $\rho$ in finite time.
\end{remark4}

\section{Blow-up}
\newtheorem {remark5}{Remark}[section]
\newtheorem{theorem5}{Theorem}[section]
\newtheorem{lemma5}{Lemma}[section]

\par
In this section, we will state two new blow-up criterions with respect to the initial data and the exact blow-up rate of strong solutions
to System (1.3).

Remark 4.1 and Lemma 4.2 imply that if we want to study the fine structure of finite
time singularities, one should assume in the following that there is a $M>0$ such that
$\|\rho(t,\cdot)\|_{L^{\infty}},\|\rho(t,\cdot)\|_{L^2}\leq e^{2Mt}\|\rho_{0}\|_{H^{s-1}}$ for all $t\in [0,T)$.
Next we will apply Lemma 4.3 to establish our first blow-up result with respect to the initial data.

\begin{theorem5} Let $z_{0}=\left(
                                                     \begin{array}{c}
                                                       u_{0} \\
                                                       \rho_{0} \\
                                                     \end{array}
                                                   \right)
\in H^s(\mathbb{R})\times H^{s-1}(\mathbb{R})$ with $s>\frac 3 2$ and
$T$ be the maximal existence time of the solution $z=\left(
                                                     \begin{array}{c}
                                                       u \\
                                                      \rho \\
                                                     \end{array}
                                                   \right) $ to System (1.3), which is guaranteed by Remark 3.1 (1).
Assume that there is a $M>0$ such that
$\|\rho(t,\cdot)\|_{L^{\infty}},\|\rho(t,\cdot)\|_{L^2}\leq e^{2Mt}\|\rho_{0}\|_{H^{s-1}}$ for all $t\in [0,T)$. Let $\varepsilon>0$ and $T^{\star}\triangleq \frac{\ln (1+\frac{2}{\varepsilon})}{\sqrt{|c|+1}(\|\rho_{0}\|_{H^{s-1}}+1)}>0$.
If there is a point $x_0\in \mathbb{R}$ such that
\begin{eqnarray*}
u'_0(x_0)<-(1+\varepsilon) K(T^{\star}),
\end{eqnarray*}
where $K(T^{\star})\triangleq \big(\frac{3(|c|+1)}{4} (e^{2MT^{\star}}(\|\rho_{0}\|_{H^{s-1}}+1))^{2}+\frac{3}{2}\, J^2(T^{\star})\big)^{\frac 1 2}>0$,
then $T<T^{\star}$. In other words, the corresponding solution to System (1.3) blows up in finite time.
\end{theorem5}

\begin{proof}
As mentioned earlier, we may assume $s=3$ here.
By (4.21) and (4.22), in view of $p\ast(\frac{3}{2}u^2)\geq 0$ and
$ \|p\ast \rho^2\|_{L^{\infty}}\leq
\|p\|_{L^{\infty}}\|\rho\|^2_{L^2}\leq
\frac{1}{2}(e^{2Mt}\|\rho_{0}\|_{H^{s-1}})^2$,
we obtain
\begin{eqnarray}
&&\ \ \ \
\frac{du_x(t,q(t,x))}{dt}\\
\nonumber&\leq&-u^2_{x}(t,q(t,x))-p\ast(\frac{c}{2}\rho^2(t,q(t,x)))+\frac{3}{2}u^2(t,q(t,x))+\frac{c}{2}\rho^2(t,q(t,x))\\
\nonumber&\leq&-u^2_{x}(t,q(t,x))+\frac{3|c|}{4} (e^{2Mt}\|\rho_{0}\|_{H^{s-1}})^{2}+\frac{3}{2}\|u(t)\|_{L^{\infty}}^2.
\end{eqnarray}
Set $m(t)\triangleq u_x (t,q(t,x_0))$ and fix $\varepsilon>0$. From (5.1) and (4.4), we have
\begin{eqnarray}
\frac{d m(t)}{dt}\leq -m^2(t)+K^2(T^{\star}), \ \ \ \forall t\in[0,T^{\star}]\cap [0,T).
\end{eqnarray}
Since  $m(0)<-(1+\varepsilon) K(T^{\star})<-K(T^{\star})$, it then follows that
$$m(t)\leq -K(T^{\star}),\ \ \ \forall t\in[0,T^{\star}]\cap [0,T).$$
By solving the inequality (5.2), we get
\begin{eqnarray}
\frac{m(0)+K(T^{\star})}{m(0)-K(T^{\star})}e^{2K(T^{\star})t}-1\leq\frac{2K(T^{\star})}{m(t)-K(T^{\star})}\leq 0.
\end{eqnarray}
Noting that $m(0)<-(1+\varepsilon) K(T^{\star})$ and $2 K(T^{\star})T^{\star}\geq \ln(1+\frac{2}{\epsilon})$, we deduce that
\begin{eqnarray}
\ln\frac{m(0)-K(T^{\star})}{m(0)+K(T^{\star})}\leq 2 K(T^{\star})T^{\star}.
\end{eqnarray}
By (5.3), (5.4) and the fact $0<\frac{m(0)+K(T^{\star})}{m(0)-K(T^{\star})}<1$, there exists
$$ 0<\,\ T\,\ <\frac{1}{2K(T^{\star})}\ln\frac{m(0)-K(T^{\star})}{m(0)+K(T^{\star})}\leq T^{\star},$$
such that $\lim\limits_{t\to T} m(t)=-\infty$. This completes the proof of the theorem.
\end{proof}

In order to establish the second blow-up result, we need the following useful lemma.
\begin{lemma5}
\cite{C-E} Let $T>0$ and $u\in C^1([0,T); H^2)$. Then for every $t\in[0,T)$, there exists at least one point $\xi(t)\in\mathbb{R}$ with
$$m(t)\triangleq \inf\limits_{x\in\mathbb{R}} (u_x(t,x))=u_x(t,\xi(t)).$$
The function $m(t)$ is absolutely continuous on $(0,T)$ with $$\frac{d m}{d t}=u_{tx}(t,\xi(t))\quad a.e.\,\ on\,\ (0,T).$$
\end{lemma5}
\begin{theorem5}
Let $z_{0}=\left(
                                                     \begin{array}{c}
                                                       u_{0} \\
                                                       \rho_{0} \\
                                                     \end{array}
                                                   \right)
\in H^s(\mathbb{R})\times H^{s-1}(\mathbb{R})$ with $s>\frac 3 2$ and
$T$ be the maximal existence time of the solution $z=\left(
                                                     \begin{array}{c}
                                                       u \\
                                                       \rho \\
                                                     \end{array}
                                                   \right) $ to System (1.3), which is guaranteed by Remark 3.1 (1).
Assume that $c\geq 0$ and the initial data satisfies that $u_0$ is odd, $\rho_0$ is even, $u^{'}_0(0)<0$ and $\rho_0(0)=0$. Then $T\leq -\frac{1}{u^{'}_0(0)}\triangleq T_0$ and $\lim\limits_{t\to T_0} u_x(t,0)=-\infty$. Moreover, if there is some $x_0\in\mathbb{R}$ such that $u_0^{'}(x_0)=\inf\limits_{x\in\mathbb{R}}u_0^{'}(x)$ and $\rho_0^{'}(x_0)\neq 0$, then there exists a $T_1\in(0,-\frac{1}{u^{'}_0(0)}]$, such that
$\limsup\limits_{t\to T_1}\big(\sup\limits_{x\in\mathbb{R}} \rho_x(t,x)\big)=+\infty$, if $\rho_0^{'}(x_0)>0$ and $\liminf\limits_{t\to T_1}\big(\inf\limits_{x\in\mathbb{R}} \rho_x(t,x)\big)=-\infty$ otherwise.
\end{theorem5}

\begin{proof}
We may assume $s=3$ here. By the assumption $u_0$ is odd, $\rho_0$ is even, and the structure of System (1.3), we have $u(t,x)$ is odd and $\rho(t,x)$ is even with respect to $x$ for $t\in(0,T)$. Thus, $u(t,0)=0$ and $\rho_x (t,0)=0$. \\
Since $\rho_0(0)=0$ and the second equation of System (1.3), it follows that
\begin{eqnarray*}
\rho(t,0)= \rho_0(0) e^{-2\int_0^t u_x(s,0)d s}=0.
\end{eqnarray*}
Set $M(t)\triangleq u_x(t,0)$. By (4.21), (4.22) and  in view of $c\geq 0$, $p\ast {u^2} \geq 0$ and
$p\ast {\rho}^2 \geq 0$, we have
\begin{eqnarray}
\frac{d M(t)}{dt}
&=&-M^2(t)-p\ast(\frac{3}{2}u^2(t,0)+\frac{c}{2}\rho^2(t,0))\\
&\leq&-M^2(t)\nonumber.
\end{eqnarray}
Note that if $M(0)=u^{'}_0(0)<0$, then $M(t)\leq M(0)<0$ for all $t\in(0,T]$.
From (5.5), we obtain $T\leq -\frac{1}{u^{'}_0(0)}$ and
\begin{eqnarray}
u_x(t,0)=M(t)\leq \frac{u^{'}_0(0)}{1+t u^{'}_0(0)} \to -\infty,
\end{eqnarray}
as $t\to -\frac{1}{u^{'}_0(0)}$.

On the other hand, applying Eq.(4.1) and differentiating the second equation in System (1.3)
with respect to $x$, we get
\begin{eqnarray}
\frac{d\rho_x(t,q(t,x))}{dt}=(-3u_x\rho_x-2u_{xx}\rho)(t,q(t,x)).
\end{eqnarray}
By Lemma 5.1, there exists $\xi(t)\in \mathbb{R}$ such that
\begin{eqnarray}
u_x(t,\xi(t))=\inf\limits_{x\in\mathbb{R}}(u_x(t,x))\quad \forall\, t\in[0,T).
\end{eqnarray}
Hence,
\begin{eqnarray}
u_{xx}(t,\xi(t))=0\quad a.e.\,\ t\in[0,T).
\end{eqnarray}
By (5.7), (5.9) and Lemma 4.1, we have
\begin{eqnarray*}
\frac{d\rho_x(t,\xi(t))}{dt}=-3u_x(t,\xi(t))\, \rho_x(t,\xi(t)),
\end{eqnarray*}
together with the assumption $u_0^{'}(x_0)=\inf\limits_{x\in\mathbb{R}}u_0^{'}(x)$ and (5.8) yields $\xi(0)=x_0$,
\begin{eqnarray}
\ \ \ \rho_x(t,\xi(t))= \rho^{'}_0(x_0) e^{-3\int_0^t u_x(s,\xi(s))d s}=\rho^{'}_0(x_0) e^{-3\int_0^t \inf\limits_{x\in\mathbb{R}}u_x(s,x)d s}.
\end{eqnarray}
Thanks to (5.6) again, we have for all $t\in [0,T)$,
$$e^{-3\int_0^t \inf\limits_{x\in\mathbb{R}}u_x(s,x)d s}\geq e^{-3\int_0^t \frac{u^{'}_0(0)}{1+s u^{'}_0(0)}d s} =\frac{1}{(1+u^{'}_0(0)t)^3}\to +\infty,$$
as $t\to -\frac{1}{u^{'}_0(0)}$.
This implies the desired result and we have completed the proof of the theorem.
\end{proof}

We conclude this section with the exact blow-up rate for blowing-up solutions to System (1.3).
\begin{theorem5}
Let $z_{0}\triangleq\left(
                                                     \begin{array}{c}
                                                       u_{0} \\
                                                       \rho_{0} \\
                                                     \end{array}
                                                   \right)
\in H^s(\mathbb{R})\times H^{s-1}(\mathbb{R})$ with $s>\frac 3 2$ and
$T<\infty$ be the blow-up time of the corresponding solution $z\triangleq\left(
                                                     \begin{array}{c}
                                                       u \\
                                                      \rho \\
                                                     \end{array}
                                                   \right) $ to System (1.3).
Assume that there is a $M>0$ such that
$\|\rho(t,\cdot)\|_{L^{\infty}},\|\rho(t,\cdot)\|_{L^2}\leq e^{2Mt}\|\rho_{0}\|_{H^{s-1}}$ for all $t\in [0,T)$. Then
\begin{eqnarray}
\lim\limits_{t\to T} \big(\inf\limits_{x\in\mathbb{R}}\{u_x(t,x)\}(T-t)\big)=-1.
\end{eqnarray}
\end{theorem5}

\begin{proof}
We may assume $s=3$ here. By the assumptions of the theorem and (4.4), we can find a $M_0>0$, such that
$$\|\rho(t,\cdot)\|_{L^\infty},\,\ \|u(t,\cdot)\|_{L^\infty}\leq M_0,\quad  \forall\, t\in [0,T).$$
Hence,
$$ \|p\ast \rho^2\|_{L^{\infty}}\leq
\|p\|_{L^1} \|\rho\|^2_{L^{\infty}}\leq {M_0}^2 \,\ and\,\
\|p\ast u^2\|_{L^{\infty}}\leq
\|p\|_{L^1} \|u\|^2_{L^{\infty}}\leq {M_0}^2,$$
which together with (4.21) and (4.22), implies that
\begin{eqnarray}
\big|\frac{d m(t)}{dt}+m^2(t)\big| \leq \widetilde{K},
\end{eqnarray}
where $\widetilde{K}=\widetilde{K}(|c|,M_0)$ is a positive constant.\\
For every $\varepsilon\in(0,\frac{1}{2})$, using the fact $m(t)\leq u_x(t,0)$ and (5.6), we can find a $t_0\in(0,T)$
such that $$m(t_0)<-\sqrt{\widetilde{K}+\frac{\widetilde{K}}{\varepsilon}}<-\sqrt{\widetilde{K}}.$$
Thanks to (5.12) again, we have $$m(t)\leq -\sqrt{\widetilde{K}}.$$
This implies that $m(t)$ is decreasing on $[t_0,T)$, hence, $$m(t)<-\sqrt{\widetilde{K}+\frac{\widetilde{K}}{\varepsilon}}<-\sqrt{\frac{\widetilde{K}}{\varepsilon}},\quad \forall\,t\in[t_0,T).$$
Noting that $-m^2(t)-\widetilde{K}\leq \frac{d m(t)}{d t}\leq -m^2(t)+\widetilde{K}\quad a.e.\,\ t\in (t_0,T),$ we get
\begin{eqnarray}
{1-\varepsilon}\leq \frac{d}{d t}\big(\frac{1}{m(t)}\big)\leq {1+\varepsilon} \quad a.e.\,\ t\in (t_0,T).
\end{eqnarray}
Integrating (5.13) with respect to $t\in [t_0,T)$ on $(t,T)$ and applying $\lim\limits_{t\to T} m(t)=-\infty$ again, we deduce that
\begin{eqnarray}
({1-\varepsilon})(T-t)\leq -\frac{1}{m(t)}\leq ({1+\varepsilon})(T-t).
\end{eqnarray}
Since $\varepsilon\in(0,\frac{1}{2})$ is arbitrary, it then follows from (5.14) that (5.11) holds.
This completes the proof of the theorem.
\end{proof}

\bigskip
\noindent\textbf{Acknowledgments} \quad  This work was partially supported by NNSFC (NO. 10971235) and the key project of Sun Yat-sen University.

\end{document}